\documentclass[reqno]{amsart}

\usepackage{amssymb,amsxtra,latexsym}

\usepackage{concrete}
\usepackage{eucal}
\usepackage[T1]{fontenc}
\usepackage{textcomp}

\usepackage[cmtip,matrix,arrow,curve]{xy}
\CompileMatrices
\newdir{((}{{\lhook\kern-.1em}}
\newdir{))}{{\rhook\kern-.1em}}
\newdir{ >}{{}*!/-5pt/@{>}}

\usepackage{booktabs}

\renewcommand\dots{\relax\ifmmode\ldots\else$\,\ldots\,$\fi}

\newcommand\note[1]%
{$^\dagger$\marginpar{\footnotesize{$^\dagger${#1}}}}

\divide\count253 by 60 
\divide\count255 by 60
\multiply\count255 by 60 
\advance\count254 by -\count255 

\def\today{\number\year-\ifnum\month<10
0\fi\number\month-\ifnum\day<10 0\fi\number\day}

\def\hour{\ifnum\count253<10
0\number\count253\else\number\count253\fi}

\def\minute{\ifnum\count254<10
0\number\count254\else\number\count254\fi}

\swapnumbers
\newtheorem{theorem}{Theorem}[section]
\newtheorem{lemma}[theorem]{Lemma}
\newtheorem{proposition}[theorem]{Proposition}
\newtheorem{corollary}[theorem]{Corollary}

\theoremstyle{definition}

\newtheorem{example}[theorem]{Example}

\newtheorem{remark}[theorem]{Remark}

\newcommand\lie{\mathfrak}

\newcommand{\g}{\lie{g}}

\newcommand{\p}{\lie{p}}
\renewcommand{\t}{\lie{t}}

\newcommand\bb[1]{{\text{\bf#1}}}

\newcommand\Z{\bb{Z}} 
\newcommand\Q{\bb{Q}}
\newcommand\R{\bb{R}} 
\newcommand\C{\bb{C}}
\renewcommand\k{\bb{k}}
\newcommand\K{\bb{K}}
\renewcommand\H{\bb{H}}
\renewcommand\P{\bb{P}}

\newcommand\ca{\mathcal}

\newcommand\X{\ca{X}}

\newcommand\func[1]{\operatorname{\mathrm{#1}}}

\newcommand\Ad{\func{Ad}}

\newcommand\End{\func{End}}

\newcommand\Hom{\func{Hom}}
\newcommand\id{\func{id}}

\newcommand\pr{\func{pr}}

\newcommand\Tor{\func{Tor}}

\newcommand\group[1]{{\text{\bf#1}}}

\newcommand\Spin{\group{Spin}}

\newcommand\SO{\group{SO}}

\newcommand\SU{\group{SU}}
\newcommand\U{\group{U}}

\newcommand\PSU{\group{PSU}}
\newcommand\A{\group{A}}
\newcommand\B{\group{B}}

\newcommand\D{\group{D}}
\newcommand\E{\group{E}}
\newcommand\F{\group{F}}
\newcommand\G{\group{G}}

\newcommand\abs[1]{\lvert#1\rvert}

\newcommand\inner[1]{\langle#1\rangle}

\newcommand\quot[1][\kern.3ex]{/\kern-.7ex/_{\kern-.4ex#1}}
\newcommand\bigquot[1][\,\,]{\big/\kern-.85ex\big/_{\!\!#1}}

\newcommand\powl{[\kern-.3ex[}
\newcommand\powr{]\kern-.3ex]}
\newcommand\bigpowl{\bigl[\kern-.6ex\bigl[}
\newcommand\bigpowr{\bigr]\kern-.6ex\bigr]}

\newcommand\mult{\mathop{\cdot}}

\newcommand\inj{\hookrightarrow}
\newcommand\sur{\mathrel{\to\kern-1.8ex\to}}
\newcommand\iso{\mathrel{\hookrightarrow\kern-1.8ex\to}}

\newcommand\longto{\longrightarrow}

\newcommand\longsur{\mathrel{\longrightarrow\kern-1.8ex\to}}
\newcommand\longiso{\mathrel{\overset{\cong}{\longrightarrow}}}

\newcommand\zerodots%
{\smash{\makebox[0pt]{$_{\lower8pt\hbox{\Large$0$}}%
\ddots^{\lower3pt\hbox{\Large$0$}}$}}}
\newcommand\bigzerodots%
{\smash{\makebox[0pt]{$_{\lower6pt\hbox{\Large$0$}}%
\ddots^{\hbox{\Large$0$}}$}}}

\newcommand\eps{\varepsilon}

\newcommand\pt{{\rm pt}}

\renewcommand\subset{\subseteq}
\renewcommand\supset{\supseteq}

\begin{document}


\title[Torsion and abelianization]{Torsion and abelianization in
equivariant cohomology}

\author{Tara S. Holm}

\email{tsh@math.cornell.edu}

\author{Reyer Sjamaar}

\email{sjamaar@math.cornell.edu}

\address{Department of Mathematics, Cornell University, Ithaca, NY
14853-4201, USA}

\thanks{The first author was partially supported by National Science
Foundation Grant DMS-0604807.  The second author was partially
supported by National Science Foundation Grant DMS-0504641.}

\dedicatory{Dedicated to Bertram Kostant on the occasion of his
eightieth birthday}

\subjclass[2000]{55N91 (14M15, 55R40, 57T15)}

\keywords{Equivariant cohomology, Schubert polynomials, homogeneous
spaces}

\date{30 March 2008}


\begin{abstract}
Let $X$ be a topological space upon which a compact connected Lie
group $G$ acts.  It is well-known that the equivariant cohomology
$H_G^*(X;\Q)$ is isomorphic to the subalgebra of Weyl group invariants
of the equivariant cohomology $H_T^*(X;\Q)$, where $T$ is a maximal
torus of $G$.  This relationship breaks down for coefficient rings
$\k$ other than $\Q$.  Instead, we prove that under a mild condition
on $\k$ the algebra $H_G^*(X,\k)$ is isomorphic to the subalgebra of
$H_T^*(X,\k)$ annihilated by the divided difference operators.
\end{abstract}


\maketitle

\thispagestyle{empty}

\tableofcontents


\section*{Introduction}\label{section;introduction}

Consider a topological space $X$ and a compact connected Lie group $G$
acting continuously on $X$.  A useful invariant of the $G$-space $X$
is Borel's equivariant cohomology algebra $H_G^*(X;\k)$ with
coefficients in a commutative ring $\k$.  Let $T$ be a maximal torus
of $G$.  It is well-known, and easy to see, how the $G$-equivariant
cohomology relates to the $T$-equivariant cohomology, at least when
$\k$ is the field of the rationals $\Q$, namely the projection
$p_X\colon X_T\to X_G$ induces an isomorphism from $H_G^*(X;\Q)$ to
the Weyl group invariants in $H_T^*(X;\Q)$.  (See for example the
standard references \cite{allday-puppe;cohomological-methods,%
brion;equivariant-cohomology-intersection-theory,%
hsiang;cohomology-theory;;1975}.)

The aim of this paper is to clarify the situation for more general
coefficient rings $\k$.  In general, the map $p_X^*$ from
$H_G^*(X;\k)$ to $H_T^*(X;\k)^W$ is neither injective nor surjective.
(See Example \ref{example;feshbach} and Section
\ref{section;invariants} for examples.)  Our first main result,
Theorem \ref{theorem;d-module}, says that the action of the Weyl group
$W$ on $H_T^*(X;\k)$ extends to an action of Demazure's algebra of
divided difference operators (which is called the nil Hecke ring in
\cite{kostant-kumar;nil-hecke;advances}).\footnote{Note added in
proof: it has recently come to our attention that this result was
obtained previously by Dale Peterson (unpublished).}  Our second main
result, Theorem \ref{theorem;d}, is valid under the assumption that a
few special primes associated with the group $G$ are invertible in
$\k$.  This assumption ensures that $p_X^*$ is injective.  Theorem
\ref{theorem;d} states that the image of $p_X^*$ consists of the
$T$-equivariant classes that are annihilated by all divided difference
operators.  However, there is a fairly natural condition, stated in
Theorem \ref{theorem;abelian}, under which the image of $p_X^*$ is
equal to $H_T^*(X;\k)^W$, namely that the roots of $G$ should not be
zero divisors in the module $H_T^*(X;\k)$.

The main applications given here are the results of Section
\ref{section;homogeneous}, which generalize classical theorems of
Borel on the cohomology of equal-rank homogeneous spaces.  In future
work we hope to apply our results to the theory of Hamiltonian
actions.

All these facts have parallels in certain other equivariant cohomology
theories, such as equivariant Chow theory and equivariant K-theory,
which we intend to explore in a sequel to this paper.  Indeed, the
Chow theory analogue of Theorem \ref{theorem;d-module} was established
by Brion in \cite{brion;equivariant-chow-torus}.

This work is based on Grothendieck's notion of the torsion index,
introduced in \cite{grothendieck;torsion-homologique}, and on the
calculus of the divided difference operators, which was initiated by
Newton and a modern version of which can be found in Demazure's papers
\cite{demazure;invariants-symetriques-entiers,%
demazure;desingularisation} and Bernstein, Gelfand and Gelfand's paper
\cite{bernstein-gelfand-gelfand;schubert}.  Our construction of the
divided difference operators as ``push-pull'' operators in
$T$-equivariant cohomology generalizes observations made by various
authors, including Aky{\i}ld{\i}z and Carrell in
\cite{akyildiz-carrell;zeros-holomorphic}, Arabia in
\cite{arabia;cycles-schubert-equivariante;inventiones}, Bressler and
Evens in \cite{bressler-evens;schubert-braid-generalized}, and Kostant
and Kumar in \cite{kostant-kumar;nil-hecke;advances}.  In fact, our
paper can be regarded as an effort to extend some of the results of
\cite{bressler-evens;schubert-braid-generalized,%
demazure;invariants-symetriques-entiers} from a space consisting of a
single point to arbitrary $G$-spaces.  We learned about the torsion
index from Totaro's recent work
\cite{totaro;torsion-index-8-other,totaro;torsion-index-spin}.  Though
we do not use much of the information on the torsion index unearthed
by him (it suffices for us to know that it exists and what its prime
factors are), we exploit some of his techniques, especially in our
results on equal-rank homogeneous spaces.

We are grateful to the American Institute of Mathematics for its
hospitality and to Joseph Wolf for several helpful discussions.  We
thank Michel Brion and Sam Evens for their incisive comments on an
earlier draft of this paper.

\section{Divided difference operators in equivariant cohomology}
\label{section;divide}

\subsection*{Equivariant cohomology}

Throughout this paper, $G$ will denote a compact connected Lie group
with maximal torus $T$ and Weyl group $W=N_G(T)/T$, $X$ will denote a
topological space equipped with a continuous (left) $G$-action, and
$\k$ will denote a commutative ring with identity.  Borel and others
have found close relationships between torsion in the integral
cohomology of $G$ and that of its classifying space $BG$.  The results
of this section and the next are applications of their work.  Many
special cases of these results are well-known and go back to Borel's
foundational paper \cite{borel;cohomologie-espaces-fibres}.  Our
presentation relies on Demazure's work
\cite{demazure;invariants-symetriques-entiers} and Bressler and Evens'
work \cite{bressler-evens;schubert-braid-generalized}, the methods of
which enable us to express the $G$-equivariant cohomology of $X$ as
the submodule of the $T$-equivariant cohomology annihilated by a
certain ideal of operators.

Let $X_G$ be the Borel homotopy quotient $EG\times^GX$ of $X$ and let
$H_G^*(X;\k)=H^*(X_G;\k)$ be the $G$-equivariant singular cohomology
with coefficients in $\k$.  Similarly, we have the $T$-equivariant
cohomology $H_T^*(X;\k)=H^*(X_T;\k)$, where $X_T=EG\times^TX$.  The
canonical map
$$
p_X\colon X_T\longsur X_G
$$
is a locally trivial fibre bundle with fibre $G/T$.  It will be
convenient to choose a preferred fibre of $p_X$.  Choose a basepoint
$e_0\in EG$ and (if $X$ is nonempty) a basepoint $x_0\in X$, and
define $i_X\colon G/T\inj X_T$ by $i_X(gT)=T(g^{-1}e_0,g^{-1}x_0)$.
Then $i_X$ is an embedding and its image is the fibre of $p_X$ over
$G(e_0,x_0)$.  If $X$ is a single point, in which case $X_T=BT$ and
$X_G=BG$, we write $i_X=i$ and $p_X=p$.  Writing $\pr_T\colon X_T\to
BT$ and $\pr_G\colon X_G\to BG$ for the canonical projections, we have
a pullback diagram
\begin{equation}\label{equation;bundles}
\vcenter{\xymatrix{
X_T\ar[r]^{p_X}\ar[d]_{\pr_T}&X_G\ar[d]^{\pr_G}\\
BT\ar[r]^p&BG,
}}
\end{equation}
in which the vertical maps are fibre bundles with fibre $X$ and the
horizontal maps are fibre bundles with fibre $G/T$.

Let $\X(T)=\Hom(T,\U(1))$ be the character group of $T$.  Let
$S=S(\X(T))$ be the symmetric algebra (over the integers) of $\X(T)$
and $S^W$ the subalgebra of Weyl group invariants.  The inclusion map
$S^W\to S$ is split injective and hence the induced map
$$S^W\otimes_\Z\k\longto S\otimes_\Z\k$$
is split injective.  For any $\Z$-module $M$, let us write
$M\otimes_\Z\k=M_\k$.  We shall identify $(S^W)_\k$ with its image in
$S_\k$.  For each character $\lambda\in\X(T)$, let $L(\lambda)$ be a
copy of $\C$ equipped with the $T$-action $(t,z)\mapsto\lambda(t)z$.
We view $L(\lambda)$ as a $T$-equivariant complex line bundle over a
point and let $L_X(\lambda)=(ET\times X\times L(\lambda))/T$ be its
extension to $X_T$.  The map $\X(T)\to H_T^2(X;\Z)$ which sends
$\lambda$ to $c_1(L_X(\lambda);\Z)$, the Chern class of
$L_X(\lambda)$, induces a ring homomorphism $S\to H_T^*(X;\Z)$.
Extending scalars to $\k$ gives a degree-doubling homomorphism of
$\k$-algebras
$$
c_X\colon S_\k\longto H_T^*(X;\k),
$$
called the \emph{characteristic homomorphism} of $X$.  The
characteristic homomorphism of a point $S_\k\to
H_T^*(\pt;\k)=H^*(BT;\k)$ is an isomorphism and from now on we shall
identify $S_\k$ with $H^*(BT;\k)$ via this isomorphism.  When no
confusion can arise, we shall write $ua$ for a product $c_X(u)a$ of
elements $u\in S_\k$ and $a\in H_T^*(X;\k)$.

\subsection*{Divided differences}

The natural Weyl group action on $X_T$ induces a $W$-action on
$H_T^*(X;\k)$ and this combines with the $S_\k$-module structure to
make $H_T^*(X;\k)$ a module over the \emph{crossed product} or
\emph{twisted group algebra} $S_\k[W]$, which is the free
$S_\k$-module on $W$ with multiplication given by
$$
\biggl(\sum_{w\in W}u_ww\biggr)\mult\biggl(\sum_{w'\in
W}u'_{w'}w'\biggr)=\sum_{w,w'\in W}u_ww(u'_{w'})ww'.
$$
Our first main result, Theorem \ref{theorem;d-module} below, says that
the $S_\k[W]$-action on $H_T^*(X;\k)$ extends naturally to an action
of the algebra of divided difference operators.  For type $A$ root
systems, these operators go back to Newton's Principia \cite[Book III,
Lemma V]{newton;principia}.  As far as we know, the definition for
general root systems goes back to
\cite{bernstein-gelfand-gelfand;schubert} and
\cite{demazure;invariants-symetriques-entiers,%
demazure;desingularisation}.  We will mostly follow the treatment
given in \cite{demazure;invariants-symetriques-entiers,%
demazure;desingularisation}, although we adopt some of the conventions
and notations of \cite{macdonald;notes-schubert-polynomials}.  Let
$R\subset\X(T)$ be the root system of $(G,T)$.  For every $\alpha\in
R$, let $G_\alpha=Z_G(\ker(\alpha))$ be the centralizer in $G$ of the
codimension one closed subgroup of $T$ fixed by $\alpha$.  Then
$G_\alpha$ has root system $\{\alpha,-\alpha\}$.  Let
$$p_{\alpha,X}\colon X_T\to X_{G_\alpha}$$
be the corresponding projection with fibre $G_\alpha/T$.  Let
$\g_\alpha$ be the Lie algebra of $G_\alpha$ and let
$\g_\C^\alpha\subset\g_\C$ be the root space of $\alpha$.  The map
$\g_\alpha\to\g_\C^\alpha$ obtained by composing the inclusion
$\g_\alpha\inj\g_\C$ with the canonical projection
$\g_\C\sur\g_\C^\alpha$ induces an $\R$-linear isomorphism
$\g_\alpha/\t\cong\g_\C^\alpha$.  This gives us a complex structure on
the tangent bundle of $G_\alpha/T\cong\P_\C^1$ and hence a natural
orientation on the fibres of $p_{\alpha,X}$.  Similarly, we orient the
fibres of $p_X$ by choosing a basis of $R$ and by identifying the
tangent space at the identity coset of $G/T$ with the sum
$\bigoplus_{\alpha\in R_+}\g_\C^\alpha$, where $R_+$ is the set of
positive roots.  (This sign convention differs from that of
\cite{bressler-evens;schubert-braid-generalized}.)  With the
orientations so defined, let
$$
p_{X,*}\colon H_T^*(X;\k)\to H_G^*(X;\k),\quad p_{\alpha,X,*}\colon
H_T^*(X;\k)\to H_{G_\alpha}^*(X;\k)
$$
be the corresponding Gysin homomorphisms.  Let
$N=\abs{R_+}=\frac12\dim(G/T)$.  Define operations $\delta$ of degree
$-2N$ and $\delta_\alpha$ of degree $-2$ on $H_T^*(X;\k)$ by
$$
\delta=p_X^*p_{X,*}\quad\text{and}\quad
\delta_\alpha=p_{\alpha,X}^*p_{\alpha,X,*}.
$$
Occasionally we shall write $\delta=\delta_X$ and
$\delta_\alpha=\delta_{\alpha,X}$ to emphasize the dependence on $X$.
These operations are functorial in two different ways.

\begin{lemma}\label{lemma;functorial}
The operations $\delta$ and $\delta_\alpha$, for every root $\alpha$,
are contravariant with respect to $G$-equivariant continuous maps and
covariant with respect to $G$-equivariant proper oriented maps.
\end{lemma}

\begin{proof}
Let $f\colon X\to Y$ a $G$-equivariant continuous map and let
$f_T\colon X_T\to Y_T$ and $f_G\colon X_G\to Y_G$ be the maps induced
by $f$.  The square
$$
\xymatrix{
X_T\ar[r]^{p_X}\ar[d]_{f_T}&X_G\ar[d]^{f_G}\\
Y_T\ar[r]^{p_Y}&**[r]Y_G
}
$$
is Cartesian, so $p_{X,*}f_T^*=f_G^*p_{Y,*}$ by the base-change
formula.  By naturality, $f_T^*p_Y^*=p_X^*f_G^*$, so
$$
\delta_Xf_T^*=p_X^*p_{X,*}f_T^*=p_X^*f_G^*p_{Y,*}
=f_T^*p_Y^*p_{Y,*}=f_T^*\delta_Y.
$$
Similarly, $\delta_{\alpha,X}f_T^*=f_T^*\delta_{\alpha,Y}$ for all
$\alpha\in R$.  If $f$ is proper and oriented, we get
$p_{Y,*}f_{T,*}=f_{G,*}p_{X,*}$ by naturality and
$f_{T,*}p_X^*=p_Y^*f_{G,*}$ by base change, which gives
$\delta_Yf_{T,*}=f_{T,*}\delta_X$.  Similarly,
$\delta_{\alpha,Y}f_{T,*}=f_{T,*}\delta_{\alpha,X}$ for all $\alpha$.
\end{proof}

Applying this result to the constant map $X\to\pt$ we obtain
$\delta_X(u\mult1)=\delta_\pt(u)\mult1$ and
$\delta_{\alpha,X}(u\mult1)=\delta_{\alpha,\pt}(u)\mult1$ for all
polynomials $u\in S_\k$, where $1\in H_T^0(X;\k)$ is the identity
element.

\begin{lemma}\label{lemma;square}
Let $\alpha$ be a root.  Then
$\delta_\alpha(p_{\alpha,X}^*(b)a)=p_{\alpha,X}^*(b)\delta_\alpha(a)$
for all $a\in H_T^*(X;\k)$ and $b\in H_{G_\alpha}^*(X;\k)$.  Hence
$\delta_\alpha(p_{\alpha,X}^*(b))=0$ and $\delta_\alpha^2=0$.
\end{lemma}

\begin{proof}
It follows from the projection formula that
$$
\delta_\alpha(p_{\alpha,X}^*(b)a)
=p_{\alpha,X}^*p_{\alpha,X,*}(p_{\alpha,X}^*(b)a)
=p_{\alpha,X}^*(bp_{\alpha,X,*}(a))
=p_{\alpha,X}^*(b)\delta_\alpha(a).
$$
In particular, $\delta_\alpha(p_{\alpha,X}^*(b))
=p_{\alpha,X}^*(b)\delta_\alpha(1)=0$, because $\delta\alpha$ is of
degree $-2$.  The identity $\delta_\alpha^2=0$ follows immediately
from this.
\end{proof}

\begin{remark}\label{remark;free}
Let $\bar{p}_X\colon X/T\to X/G$, $q_T\colon X_T\to X/T$ and
$q_G\colon X_G\to X/G$ be the natural projections.  If $G$ acts freely
on $X$, then $X/T$ is an oriented fibre bundle over $X/G$, so we can
define an operator $\bar{\delta}$ on $H^*(X/T;\k)$ by
$\bar{\delta}=\bar{p}_X^*\bar{p}_{X,*}$.  Moreover, the square
$$
\xymatrix{
X_T\ar[r]^{p_X}\ar[d]_{q_T}&X_G\ar[d]^{q_G}\\
X/T\ar[r]^{\bar{p}_X}&X/G
}
$$
is Cartesian and $q_T^*$ is an isomorphism.  Arguing as in the proof
of Lemma \ref{lemma;functorial} we find $\delta
q_T^*=q_T^*\bar{\delta}$, so the natural transformations $\delta$ and
$\bar{\delta}$ are isomorphic.  For this reason we shall identify the
natural transformations $\delta$ and $\bar{\delta}$ on the category of
free $G$-spaces.  Likewise, for any root $\alpha$, we let
$\bar{p}_{\alpha,X}\colon X/T\to X/G_\alpha$ be the natural projection
and define
$\bar{\delta}_\alpha=\bar{p}_{\alpha,X}^*\bar{p}_{\alpha,X,*}$ for a
free $G_\alpha$-space $X$.  Then $\delta_\alpha
q_T^*=q_T^*\bar{\delta}_\alpha$.  For an arbitrary $G$-space $X$, it
is plain from the definitions that $\delta_X=\bar{\delta}_{EG\times
X}$ and $\delta_{\alpha,X}=\bar{\delta}_{\alpha,EG\times X}$.  In this
manner, all questions regarding the operations $\delta$ and
$\delta_\alpha$ can be reduced to the case of a free action.
\end{remark}

Next we show that the Weyl group action can be expressed in terms of
the operators $\delta_\alpha$ and that the $\delta_\alpha$ satisfy an
integration by parts rule.

\begin{proposition}\label{proposition;divide-leibniz}
Let $\alpha$ be a root.  Then
\begin{enumerate}
\item\label{item;divide}
$s_\alpha=1-\alpha\delta_\alpha$, where $s_\alpha\in W$ denotes the
reflection in $\alpha$ and $1$ denotes the identity operator on
$H_T^*(X;\k)$;
\item\label{item;leibniz}
$\delta_\alpha(a_1a_2)=\delta_\alpha(a_1)a_2
+s_\alpha(a_1)\delta_\alpha(a_2)$ for all $a_1$, $a_2\in H_T^*(X;\k)$.
\end{enumerate}
\end{proposition}

\begin{proof}
We may assume without loss of generality that $G=G_\alpha$.  Then
$R=\{\alpha,-\alpha\}$, $p_{\alpha,X}=p_X$ and $\delta_\alpha=\delta$.
We may also assume that $G$ acts freely on $X$.  (If it does not,
replace $X$ with $EG\times X$; see Remark \ref{remark;free}.)  We will
prove the proposition by pulling up the fibre bundle $X/T\to X/G$ along the
projection map $p_X$.  Let $X'=X\times^TG$ be the quotient of $X\times
G$ by the free $T$-action $t\mult(x,g)=(tx,tg)$.  Then $G$ acts on
$X'$ by right multiplication on the fibre $G$.  Let
$$
c_{X'}\colon S_\k\to H^*(X'/T;\k)
$$
be the characteristic homomorphism for $X'$ considered as a $T$-space.
The quotient of $X'$ by $G$ is $X'/G=(X\times^TG)/G\cong X/T$.  With
this identification, the map $p_{X'}\colon X'/T\to X'/G$ becomes the
map $f\colon X'/T\to X/T$ given by $f([x,g])=[x]$.  Here $[x,g]\in
X'/T$ denotes the $T$-orbit of $(x,g)\in X\times G$ and $[x]\in X/T$
denotes the $T$-orbit of $x\in X$.  We define a second map $h\colon
X'/T\to X/T$ by $h([x,g])=[g^{-1}x]$.  The square
\begin{equation}\label{equation;change-bundle}
\vcenter{\xymatrix{
X'/T\ar[r]_h\ar[d]^f&X/T\ar[d]^{p_X}\ar@/_/@<-0.5ex>@{.>}[l]_\tau\\
X/T\ar[r]^{p_X}\ar@/^/@<0.5ex>@{.>}[u]^\sigma&X/G
}}
\end{equation}
is Cartesian.  The section $\tau$ of $h$ is defined by
$\tau([x])=[x,1]$ and the section $\sigma$ of $f$ is defined by
$\sigma(x)=[x,n_\alpha]$, where $n_\alpha\in N_G(T)$ is a
representative of $s_\alpha$.  Lemma \ref{lemma;functorial} gives
\begin{equation}\label{equation;push-pull-upstairs}
h^*\delta_X=\delta_{X'}h^*,
\end{equation}
where $\delta_{X'}=p_{X'}^*p_{X',*}=f^*f_*$.  Note that $h$ is induced
by the $G$-equivariant map $X'\to X$ defined by $[x,g]\mapsto
g^{-1}x$.  It follows that $h$ is equivariant with respect to the
given Weyl group action on $X/T$ and the fibrewise Weyl group action
on $X'/T$ defined by $w\mult[x,g]=[x,gn^{-1}]$, where $n\in N_G(T)$ is
a representative of $w\in W$.  Therefore
\begin{equation}\label{equation;reflect-upstairs}
h^*(s_\alpha(a))=s_\alpha(h^*(a))
\end{equation}
for all $a\in H^*(X/T;\k)$.  Because $h$ is induced by an equivariant
map, we have also
\begin{equation}\label{equation;pull-characteristic}
h^*(c_X(u)a)=c_{X'}(u)h^*(a)
\end{equation}
for all $u\in S_\k$ and $a\in H^*(X/T;\k)$.  Since $h$ has a section,
$h^*$ is injective, so because of
\eqref{equation;push-pull-upstairs}--\eqref{equation;pull-characteristic}
the proposition will follow from
\begin{equation}\label{equation;divide-leibniz-upstairs}
s_\alpha(b)=b-c_{X'}(\alpha)\delta_{X'}(b),\qquad
\delta_{X'}(b_1b_2)=\delta_{X'}(b_1)b_2+s_\alpha(b_1)\delta_{X'}(b_2)
\end{equation}
for all $b$, $b_1$, $b_2\in H^*(X'/T;\k)$.  Let $\xi=\sigma_*(1)\in
H^2(X'/T;\k)$.  See \cite[Sections
III.3--5]{bott-samelson;applications-theory-morse-symmetric} or
\cite[Section 2]{demazure;desingularisation} for the following facts:
\begin{equation}\label{equation;bundle-relation}
\begin{split}
\delta_{X'}(1)=0,\quad\delta_{X'}(\xi)=1,&\quad s_\alpha(1)=1,\quad
s_\alpha(\xi)=-\xi-f^*(c_X(\alpha)),\\
\xi^2+f^*(c_X(\alpha))\xi=0,&\qquad
c_{X'}(\alpha)=f^*(c_X(\alpha))+2\xi.
\end{split}
\end{equation}
Since $f$ is $W$-invariant, we have $s_\alpha(f^*(a))=f^*(a)$ for all
$a\in H^*(X/T;\k)$.  By the Leray-Hirsch theorem, the classes $1$ and
$\xi$ form a basis of $H^*(X'/T;\k)$ regarded as an
$H^*(X/T;\k)$-module via the map $f^*$.  Let $b\in H^*(X'/T;\k)$.
Writing $b=f^*(a_1)+f^*(a_2)\xi$ with $a_1$, $a_2\in H^*(X/T;\k)$ and
using \eqref{equation;bundle-relation}, we find
\begin{align*}
s_\alpha(b)&=f^*(a_1)-f^*(a_2)(\xi+f^*(c_X(\alpha))),\\
b-c_{X'}(\alpha)\delta_{X'}(b)&=
f^*(a_1)+f^*(a_2)\xi-(f^*(c_X(\alpha))+2\xi)f^*(a_2),
\end{align*}
so $s_\alpha(b)=b-c_{X'}(\alpha)\delta_{X'}(b)$.  The second equality
in \eqref{equation;divide-leibniz-upstairs} follows from
\eqref{equation;bundle-relation} in a similar way.
\end{proof}

For any $u\in S=H_T^*(\pt;\Z)$, the difference $u-s_\alpha(u)$ is
uniquely divisible by $\alpha$, so it follows from Proposition
\ref{proposition;divide-leibniz}\eqref{item;divide} that
$\delta_{\alpha,\pt}$ is a divided difference operator on the
polynomial ring $S$.

\begin{corollary}\label{corollary;divide}
Let $\k=\Z$ and let $\pt$ denote the space consisting of one point.
Then
$$
\delta_{\alpha,\pt}(u)=\frac{u-s_\alpha(u)}\alpha
$$
for all $u\in S$.
\end{corollary}

This formula can also be found in
\cite{akyildiz-carrell;zeros-holomorphic} and
\cite{bressler-evens;schubert-braid-generalized}.  It follows from
this that, for an arbitrary coefficient ring $\k$, the operator
$\delta_{\alpha,\pt}$ on $S_\k$ is equal to
$\alpha^{-1}(1-s_\alpha)\otimes\id_\k$.

\begin{lemma}\label{lemma;reflect}
The identities
$$
s_\alpha\delta_\alpha=\delta_\alpha,\quad\delta_\alpha
s_\alpha=-\delta_\alpha,\quad\delta_{-\alpha}=-\delta_\alpha,\quad
w\delta_\alpha w^{-1}=\delta_{w(\alpha)}
$$
hold for all $\alpha\in R$ and $w\in W$.
\end{lemma}

\begin{proof}
The identity $s_\alpha\delta_\alpha=\delta_\alpha$ follows from Lemma
\ref{lemma;square} and Proposition
\ref{proposition;divide-leibniz}\eqref{item;divide}.  From Corollary
\ref{corollary;divide} we get $\delta_\alpha(\alpha)=2$, whence, by
Lemma \ref{lemma;square} and Proposition
\ref{proposition;divide-leibniz}\eqref{item;leibniz},
$$
\delta_\alpha s_\alpha=\delta_\alpha(1-\alpha\delta_\alpha)
=\delta_\alpha-2\delta_\alpha=-\delta_\alpha.
$$
Replacing $\alpha$ with $-\alpha$ has the effect of reversing the
complex structure on $G_\alpha/T$ and thus reversing the sign of
$p_{\alpha,X,*}$.  Therefore $\delta_{-\alpha}=-\delta_\alpha$.
Finally, for any Weyl group element $w$, conjugation by $w$ maps
$G_\alpha$ to $G_{w(\alpha)}$ and induces a biholomorphic map
$G_\alpha/T\to G_{w(\alpha)}/T$.  The action of $w$, viewed as a
homeomorphism $X_T\to X_T$, induces a homeomorphism $\bar{w}\colon
X_{G_\alpha}\to X_{G_{w(\alpha)}}$, and the diagram
$$
\xymatrix{
G_\alpha/T\ar[r]\ar[d]_{\Ad(w)}&
**[r]X_T\ar[r]^{p_{\alpha,X}}\ar[d]_w&X_{G_\alpha}\ar[d]_{\bar{w}}\\
G_{w(\alpha)}/T\ar[r]&X_T\ar[r]^{p_{w(\alpha),X}}&
**[r]X_{G_{w(\alpha)}}
}
$$
commutes.  The square on the right is Cartesian, so
$\bar{w}p_{w(\alpha),X,*}=p_{\alpha,X,*}w$.  By naturality,
$wp_{w(\alpha),X}^*=p_{\alpha,X}^*\bar{w}$, and so $w\delta_\alpha
w^{-1}=\delta_{w(\alpha)}$.
\end{proof}

Our next result relates the operations $\delta_\alpha$ to the
operation $\delta$.  The proof is based on a commutative diagram
\begin{equation}\label{equation;resolution-bundle}
\vcenter{\xymatrix@=1.5em{
G/T\ar[rrr]^{j'}\ar[dd]&&&X'/T\ar[rr]^h\ar[dd]^f&&X/T\ar[dd]^{p_X}\\
&\tilde{G}/T\ar[r]^{\tilde{\jmath}}\ar[ul]_q\ar[dl]
&\tilde{X}/T\ar[ur]^{q_X}\ar[dr]_{\tilde{f}}&&&\\
\pt\ar[rrr]^j&&&X/T\ar[rr]^{p_X}&&X/G,
}}
\end{equation}
which is defined for every free $G$-space $X$ as follows.  The square
on the right is as in \eqref{equation;change-bundle}.  The object
$\tilde{G}/T$ is a \emph{Bott-Samelson variety}.  We recall the
construction from
\cite[Section~III.4]{bott-samelson;applications-theory-morse-symmetric}.
Choose a reduced decomposition
\begin{equation}\label{equation;decomposition}
w_0=s_{\beta_1}s_{\beta_2}\cdots s_{\beta_N}
\end{equation}
of the longest Weyl group element $w_0$ in terms of simple roots
$\beta_1$, $\beta_2$,\dots, $\beta_N$.  Put $\alpha_1=\beta_1$ and
$$
\alpha_k=s_{\beta_1}s_{\beta_2}\cdots s_{\beta_{k-1}}(\beta_k)
$$
for $2\le k\le N$.  Then $\{\alpha_1,\alpha_2,\dots,\alpha_N\}=R_+$.
For $1\le k\le N$ define
$$
G^{(k)}
=G_{\alpha_1}\times^TG_{\alpha_2}\times^T\cdots\times^TG_{\alpha_k},
$$
where the $n$th copy of $T$ acts on $G_{\alpha_n}\times
G_{\alpha_{n+1}}$ by $t\mult(g_n,g_{n+1})=(g_nt^{-1},tg_{n+1})$.  Put
$\tilde{G}=G^{(N)}$.  The torus $T$ acts freely on $\tilde{G}$ by
right multiplication on $G_{\alpha_N}$.  Let $q\colon\tilde{G}/T\to
G/T$ be the map induced by multiplication,
$$[g_1,g_2,\dots,g_N]\longmapsto[g_1g_2\cdots g_N].$$
The map $q\colon\tilde{G}/T\to G/T$ is the \emph{Bott-Samelson
resolution} of the flag variety associated with the given reduced
decomposition of $w_0$, so named because it is a simultaneous
resolution of singularities of all the Schubert varieties in $G/T$.
(See \cite[Sections 3--4]{demazure;desingularisation}.)  The space
$\tilde{X}/T$ is obtained by replacing the fibre $G/T$ of the bundle
$f$ by the Bott-Samelson variety $\tilde{G}/T$.  Namely, let $X^{(k)}$
be the $T$-space $X\times^TG^{(k)}$ and let $\tilde{X}=X^{(N)}$.  We
define $q_X\colon\tilde{X}/T\to X'/T$ to be the map
$$[x,g_1,g_2,\dots,g_N]\longmapsto[x,g_1g_2\cdots g_N]$$
and $\tilde{f}\colon\tilde{X}/T\to X/T$ to be the composition of $q_X$
and $f$.  The map $j\colon\pt\to X/T$ is the inclusion of an
arbitrarily chosen point in $X/T$ and $j'\colon G/T\to X'/T$ and
$\tilde{\jmath}\colon\tilde{G}/T\to\tilde{X}/T$ are the inclusions of
the fibres of $f$, resp.\ $\tilde{f}$, over this point.

A useful property of the map $\tilde{f}$ is that it splits into a
sequence of $\P^1$-bundles with sections, which we shall call the
\emph{Bott-Samelson tower} of $X$:
$$
\xymatrix{
\tilde{X}/T\ar@<0.5ex>[r]^-{\pi_N}&
\cdots\ar@<0.5ex>[r]^-{\pi_{k+1}}\ar@<0.5ex>[l]^-{\sigma_N}&
X^{(k)}/T\ar@<0.5ex>[r]^-{\pi_k}\ar@<0.5ex>[l]^-{\sigma_{k+1}}&
X^{(k-1)}/T\ar@<0.5ex>[r]^-{\pi_{k-1}}\ar@<0.5ex>[l]^-{\sigma_k}&
\cdots\ar@<0.5ex>[r]^-{\pi_2}\ar@<0.5ex>[l]^-{\sigma_{k-1}}&
X^{(1)}/T\ar@<0.5ex>[r]^-{\pi_1}\ar@<0.5ex>[l]^-{\sigma_2}&
X/T\ar@<0.5ex>[l]^-{\sigma_1}.
}
$$
The projections $\pi_k$ and the sections $\sigma_k$ are defined by
\begin{equation}\label{equation;tower}
\begin{split}
\pi_k([x,g_1,g_2,\dots,g_{k-1},g_k])=[x,g_1,g_2,\dots,g_{k-1}],\\
\sigma_k([x,g_1,g_2,\dots,g_{k-1}])=[x,g_1,g_2,\dots,g_{k-1},n_k],
\end{split}
\end{equation}
where $n_k$ represents the nontrivial element $s_{\alpha_k}$ of the
Weyl group of $G_{\alpha_k}$.

\begin{proposition}\label{proposition;long}
For every reduced decomposition \eqref{equation;decomposition} of the
longest Weyl group element $w_0$, we have
$\delta=\delta_{\beta_1}\delta_{\beta_2}\cdots\delta_{\beta_N}$.
\end{proposition}

\begin{proof}
Let us assume, as we may, that $G$ acts freely on $X$.  First consider
the case where $X$ is $G$ acting on itself by left multiplication.  In
that case the assertion is that
\begin{equation}\label{equation;long-g}
\delta_G=\delta_{\beta_1,G}\delta_{\beta_2,G}\cdots\delta_{\beta_N,G}
\end{equation}
as operators on $H_T^*(G;\k)\cong H^*(G/T;\k)$.  Because $H^*(G/T;\Z)$
is a free and finitely generated abelian group, it suffices to check
this identity over the ring $\k=\Z$.  Both sides vanish except in
degree $2N$, so it is enough to show that they agree on the
orientation class $\theta\in H^{2N}(G/T;\Z)$.  Let
$$d=\prod_{\alpha\in R_+}\alpha\in S^N$$
be the discriminant of $G$.  By
\cite[\S~4]{demazure;invariants-symetriques-entiers} the image of $d$
under the characteristic homomorphism $c_G=i^*\colon S\to H^*(G/T;\Z)$
is
\begin{equation}\label{equation;discriminant}
i^*(d)=\abs{W}\theta.
\end{equation}
From Corollary \ref{corollary;divide} and \cite[Proposition
3]{demazure;invariants-symetriques-entiers} we obtain the identity
\begin{equation}\label{equation;antisymmetry}
\delta_{\beta_1,\pt}\delta_{\beta_2,\pt}\cdots\delta_{\beta_N,\pt}
=\frac1d\sum_{w\in W}\det(w)w,
\end{equation}
valid over the ring $\Z$.  Since $w(d)=\det(w)d$, we find
$$
\delta_{\beta_1,\pt}\delta_{\beta_2,\pt}\cdots\delta_{\beta_N,\pt}(d)
=\frac1d\sum_{w\in W}\det(w)w(d)=\abs{W}.
$$
By naturality and by \eqref{equation;discriminant} this gives
\begin{align*}
\abs{W}&=i^*\delta_{\beta_1,\pt}\delta_{\beta_2,\pt}\cdots
\delta_{\beta_N,\pt}(d)\\
&=\delta_{\beta_1,G}\delta_{\beta_2,G}
\cdots\delta_{\beta_N,G}i^*(d)\\
&=\abs{W}\delta_{\beta_1,G}\delta_{\beta_2,G}\cdots
\delta_{\beta_N,G}(\theta).
\end{align*}
Dividing by $\abs{W}$ yields
$$
\delta_{\beta_1,G}\delta_{\beta_2,G}\cdots\delta_{\beta_N,G}(\theta)
=1=\delta_G(\theta)
$$
and hence \eqref{equation;long-g}.

We will deduce from this that the proposition is true for an arbitrary
free $G$-space $X$.  The following square contained in the diagram
\eqref{equation;resolution-bundle} is Cartesian:
$$
\xymatrix{
\tilde{G}/T\ar[r]^{\tilde{\jmath}}\ar[d]_q&\tilde{X}/T\ar[d]^{q_X}\\
G/T\ar[r]^{j'}&X'/T.
}
$$
As in the proof of Lemma \ref{lemma;functorial}, this implies the
identities
\begin{equation}\label{equation;delta-resolution}
q_*\tilde{\jmath}^*=j^{\prime,*}q_{X,*},\qquad
\tilde{\jmath}^*\tilde{\delta}_X=\tilde{\delta}_G\tilde{\jmath}^*,
\qquad\tilde{\jmath}^*\tilde{\delta}_{\alpha,X}
=\tilde{\delta}_{\alpha,G}\tilde{\jmath}^*
\end{equation}
for all roots $\alpha$.  Here $\tilde{\delta}_X$ and
$\tilde{\delta}_{\alpha,X}$ are operators on $H^*(\tilde{X}/T;\k)$ and
$\tilde{\delta}_G$ and $\tilde{\delta}_{\alpha,G}$ are operators on
$H^*(\tilde{G}/T;\k)$ defined by
\begin{alignat*}{2}
\tilde{\delta}_X&=\tilde{f}^*\tilde{f}_*=q_X^*\delta_{X'}q_{X,*},&
\qquad\tilde{\delta}_{\alpha,X}&=q_X^*\delta_{\alpha,X'}q_{X,*},\\
\tilde{\delta}_G&=q^*\delta_Gq_*,&\qquad
\tilde{\delta}_{\alpha,G}&=q^*\delta_{\alpha,G}q_*.
\end{alignat*}
For $a\in H^*(X/T;\k)$ and $b\in H^*(\tilde{X}/T;\k)$ we have
\begin{multline*}
\tilde{\delta}_{\alpha,X}(\tilde{f}^*(a)b)
=q_X^*\delta_{\alpha,X'}q_{X,*}(q_X^*(f^*(a))b)
=q_X^*\delta_{\alpha,X'}(f^*(a)q_{X,*}(b))\\
=q_X^*(f^*(a)\delta_{\alpha,X'}(q_{X,*}(b)))
=\tilde{f}^*(a)q_X^*\delta_{\alpha,X'}q_{X,*}(b)
=\tilde{f}^*(a)\tilde{\delta}_{\alpha,X}(b),
\end{multline*}
where we used the $H^*(X/T;\k)$-linearity of $\delta_{\alpha,X'}$
(Lemma \ref{lemma;square}) and the projection formula.  The same
argument applies to $\tilde{\delta}_X$.  Thus $\tilde{\delta}_X$ and
$\tilde{\delta}_{\alpha,X}$ are linear over $H^*(X/T;\k)$:
\begin{equation}\label{equation;linear}
\tilde{\delta}_X(\tilde{f}^*(a)b)=\tilde{f}^*(a)\tilde{\delta}_X(b),
\qquad\tilde{\delta}_{\alpha,X}(\tilde{f}^*(a)b)
=\tilde{f}^*(a)\tilde{\delta}_{\alpha,X}(b)
\end{equation}
for all $a\in H^*(X/T;\k)$ and $b\in H^*(\tilde{X}/T;\k)$.

The Bott-Samelson map $q\colon\tilde{G}/T\to G/T$, being birational,
has mapping degree $1$.  This implies $q_*(1)=1$ and hence
$1=q_*\tilde{\jmath}^*(1)=j^{\prime,*}q_{X,*}(1)$ by the first formula
in \eqref{equation;delta-resolution}.  Since $j^{\prime,*}$ is
injective in degree $0$, this gives $q_{X,*}(1)=1$.  Therefore
\begin{equation}\label{equation;resolution}
q_*q^*=\id,\qquad q_{X,*}q_X^*=\id
\end{equation}
by the projection formula.  By left multiplying both sides of
\eqref{equation;long-g} by $q^*$ and right multiplying by $q_*$ and
using \eqref{equation;resolution}, we find $\tilde{\delta}_G=
\tilde{\delta}_{\beta_1,G}\tilde{\delta}_{\beta_2,G}\cdots
\tilde{\delta}_{\beta_N,G}$.  In particular,
\begin{equation}\label{equation;orientation-resolution}
1=\tilde{\delta}_G(\tilde{\theta})
=\tilde{\delta}_{\beta_1,G}\tilde{\delta}_{\beta_2,G}\cdots
\tilde{\delta}_{\beta_N,G}(\tilde{\theta}),
\end{equation}
where $\tilde{\theta}=q^*(\theta)$.  Since $q$ has mapping degree $1$,
$\tilde\theta$ is the orientation class of the variety $\tilde{G}/T$.
We assert that $\tilde{\theta}$ extends to a class in
$H^{2N}(\tilde{X}/T;\k)$ which is a product of degree $2$ classes.
These classes are defined as follows.  Let $\sigma_k$ be the section
\eqref{equation;tower} of the Bott-Samelson tower.  The class
$\sigma_{k,*}(1)$ is the Thom class of the $k$th ``storey'' of the
tower.  Let $\xi_k\in H^2(\tilde{X}/T;\k)$ be the pullback of this
class,
$$
\xi_k=(\pi_N)^*(\pi_{N-1})^*\cdots(\pi_{k+1})^*\sigma_{k,*}(1),
$$
and for each subset $K$ of $\{1,2,\dots,N\}$ let 
$$\xi_K=\prod_{k\in K}\xi_k\in H^{2\abs{K}}(\tilde{X}/T;\k).$$
It follows from the Leray-Hirsch theorem, applied to each storey of
the Bott-Samel\-son tower, that the $2^N$ classes $\xi_K$ form a basis
of $H^*(\tilde{X}/T;\k)$ viewed as an $H^*(X/T;\k)$-module via the map
$\tilde{f}^*$.  The top-degree basis element is
$$
\xi_{\{1,2,\dots,N\}}=\xi_1\xi_2\cdots\xi_N\in
H^{2N}(\tilde{X}/T;\k),
$$
which restricts to the orientation class of the fibre $\tilde{G}/T$,
i.e.\ $\tilde{\jmath}^*(\xi_1\xi_2\cdots\xi_N)=\tilde{\theta}$.
Substituting this into \eqref{equation;orientation-resolution} and
using the second formula in \eqref{equation;delta-resolution} we get
$$
\tilde{\jmath}^*\tilde{\delta}_X(\xi_1\xi_2\cdots\xi_N)
=\tilde{\jmath}^*\tilde{\delta}_{\beta_1,X}\tilde{\delta}_{\beta_2,X}
\cdots\tilde{\delta}_{\beta_N,X}(\xi_1\xi_2\cdots\xi_N).
$$
Both sides of this identity are of degree $0$, so
\begin{equation}\label{equation;orientation}
\tilde{\delta}_X(\xi_1\xi_2\cdots\xi_N)
=\tilde{\delta}_{\beta_1,X}\tilde{\delta}_{\beta_2,X}
\cdots\tilde{\delta}_{\beta_N,X}(\xi_1\xi_2\cdots\xi_N).
\end{equation}

Finally, we can finish the proof of the proposition.  Since $h^*$ is
injective (because $h$ has a section; see
\eqref{equation;change-bundle}), the desired result will follow from
$$
\delta_{X'}=\delta_{\beta_1,X'}\delta_{\beta_2,X'}\cdots
\delta_{\beta_N,X'}.
$$
It follows from \eqref{equation;resolution} that $q_X^*$ is injective,
so in fact it is enough to show that
\begin{equation}\label{equation;bott-samelson}
\tilde{\delta}_X=\tilde{\delta}_{\beta_1,X}\tilde{\delta}_{\beta_2,X}
\cdots\tilde{\delta}_{\beta_N,X}.
\end{equation}
Because of the linearity property \eqref{equation;linear}, it suffices
to check \eqref{equation;bott-samelson} on the basis elements $\xi_K$.
Since both sides of \eqref{equation;bott-samelson} are operators of
degree $-2N$, they annihilate all classes $\xi_K$ except when
$K=\{1,2,\dots,N\}$, in which case \eqref{equation;bott-samelson} is
equivalent to \eqref{equation;orientation}.  This proves
\eqref{equation;bott-samelson}.
\end{proof}

\subsection*{The operators $\partial_w$}

Let $\alpha$ and $\beta$ be simple roots and let $m$ be the order of
$s_\alpha s_\beta$.  Then $s_\alpha s_\beta\cdots=s_\beta
s_\alpha\cdots$ ($m$ factors on both sides) is the longest element of
the Weyl group of the root subsystem of $R$ generated by $\alpha$ and
$\beta$, so from Proposition \ref{proposition;long} we obtain the
\emph{braid relation}
\begin{equation}\label{equation;braid}
\delta_\alpha\delta_\beta\cdots=\delta_\beta\delta_\alpha\cdots,
\end{equation}
where each product has $m$ factors.  Let $w$ be a Weyl group element
of length $k$ and let $w=s_{\alpha_1}s_{\alpha_2}\cdots s_{\alpha_k}$
be a reduced decomposition of $w$ in terms of simple roots $\alpha_1$,
$\alpha_2$,\dots, $\alpha_k$.  Define an endomorphism $\partial_w$ of
$H_T^*(X;\k)$ of degree $-2k$ by
\begin{equation}\label{equation;partial}
\partial_w=\delta_{\alpha_1}\delta_{\alpha_2}\cdots\delta_{\alpha_k}.
\end{equation}
This operator depends on the choice of the basis of $R$, but it
follows from the braid relations and Matsumoto's theorem
(\cite[\S~IV.1, Proposition~5]{bourbaki;groupes-algebres}) that it
does not depend on the reduced decomposition of $w$.  Note that
\begin{equation}\label{equation;integrate}
\partial_{w_0,X}=\delta_X=p_X^*p_{X,*}.
\end{equation}
From \eqref{equation;antisymmetry} we get the following identity,
valid over the integers,
\begin{equation}\label{equation;long-symmetric}
\partial_{w_0,\pt}=\frac1d\sum_{w\in W}\det(w)w.
\end{equation}

As in \cite[Section~4]{demazure;invariants-symetriques-entiers}, the
braid relations \eqref{equation;braid} and the relations
$\delta_\alpha^2=0$ (Lemma \ref{lemma;square}) imply the following
relations among the $\partial_w$:
\begin{equation}\label{equation;relation}
\partial_w\partial_{w'}=\partial_{ww'}\quad\text{if
$l(ww')=l(w)+l(w')$},
\quad\partial_w\partial_{w'}=0\quad\text{otherwise}.
\end{equation}
From these relations and Proposition
\ref{proposition;divide-leibniz}\eqref{item;leibniz} one can deduce a
``Leibniz'' rule for $\partial_w$, which is stated in \cite[Chapter
II]{macdonald;notes-schubert-polynomials}.  This reference is hard to
find and it treats only the case of a one-point space and the
classical root system of type $A$, so let us review the argument.  Let
$$\mu\colon H_T^*(X;\k)\otimes_\k H_T^*(X;\k)\to H_T^*(X;\k)$$
be the multiplication map.  Then Proposition
\ref{proposition;divide-leibniz}\eqref{item;leibniz} is equivalent to
$$
\delta_\alpha\circ\mu
=\mu\circ(\delta_\alpha\otimes1+s_\alpha\otimes\delta_\alpha).
$$
Substituting \eqref{equation;partial} and writing $s_j$ for the simple
reflection $s_{\alpha_j}$ gives
$$
\partial_w\circ\mu=\mu\circ(\partial_{s_1}\otimes1
+s_1\otimes\partial_{s_1})\circ
(\partial_{s_2}\otimes1+s_2\otimes\partial_{s_2})\circ\cdots\circ
(\partial_{s_k}\otimes1+s_k\otimes\partial_{s_k}).
$$
This expands to an expression of the form
$$
\partial_w\circ\mu
=\mu\circ\sum_\bb{t}\phi(\bb{s},\bb{t})\otimes\partial_{\bb{t}}.
$$
Here $\bb{s}$ is the word $(s_1,s_2,\dots,s_k)$; the sum is over all
subwords $\bb{t}=(t_1,t_2,\dots,t_l)$ of $\bb{s}$ of $l$ letters,
where $0\le l\le k$; $\partial_{\bb{t}}$ denotes
$\partial_{t_1}\partial_{t_2}\cdots\partial_{t_l}$; and
$$
\phi(\bb{s},\bb{t})=\phi_1(\bb{s},\bb{t})\phi_2(\bb{s},\bb{t})\cdots
\phi_k(\bb{s},\bb{t})
$$
with
$$
\phi_j(\bb{s},\bb{t})=
\begin{cases}
s_j&\text{if $s_j\in\bb{t}$},\\
\partial_{s_j}&\text{if $s_j\not\in\bb{t}$},
\end{cases}
$$
for $1\le j\le k$.  If the subword $\bb{t}$ is not reduced, then
$\partial_{\bb{t}}=0$ by \eqref{equation;relation}, so we can write
$$
\partial_w\circ\mu=\mu\circ\sum_{w'\le
w}w'\partial_{w/w'}\otimes\partial_{w'}.
$$
Here $\partial_{w/w'} =(w')^{-1}\sum_{\bb{t}}\phi(\bb{s},\bb{t})$,
summed over all reduced subwords $\bb{t}$ of $\bb{s}$ such that
$w'=t_1t_2\cdots t_l$.  The operator $\partial_{w/w'}$ has degree
$l(w)-l(w')$ (viewed as an endomorphism of $S_\k$) and it follows from
the independence of the $\partial_w$ that $\partial_{w/w'}$ does not
depend on the reduced decomposition of $w$.  This proves the following
product rule.

\begin{lemma}\label{lemma;leibniz}
$\partial_w(a_1a_2)=\sum_{w'\le
w}w'(\partial_{w/w'}(a_1))\partial_{w'}(a_2)$ for all $a_1$ and $a_2$
in $H_T^*(X;\k)$.
\end{lemma}

\subsection*{Demazure's algebra}

Let $\ca{D}_\k$ be the algebra of endomorphisms of
$S_\k=H_T^*(\pt;\k)$ generated by the operators $\delta_{\alpha,\pt}$
for $\alpha\in R$ and the elements of $S_\k$ (regarded as
multiplication operators).  The generators are $(S^W)_\k$-linear, so
$\ca{D}_\k$ is a subalgebra of $\End_{(S^W)_\k}(S_\k)$.  By
Proposition \ref{proposition;divide-leibniz}\eqref{item;divide},
$\ca{D}_\k$ contains the Weyl group and hence the twisted group
algebra $S_\k[W]$.

\begin{proposition}\label{proposition;d}
The algebra $\ca{D}_\k$ is the free left $S_\k$-module on the elements
$(\partial_w)_{w\in W}$ equipped with the multiplication rule
$$
\biggl(\sum_{w\in W}u_w\partial_w\biggr)\mult\biggl(\sum_{w'\in
W}u'_{w'}\partial_{w'}\biggr)=\sum_{w,w'\in W}
\sum_{w''}u_ww''\bigl(\partial_{w/w''}(u'_{w'})\bigr)\partial_{w''w'},
$$
where the second sum on the right ranges over all $w''\in W$
satisfying $w''\le w$ and $l(w''w')=l(w'')+l(w')$.  
\end{proposition}

\begin{proof}
That the $\partial_w$ are a basis of $\ca{D}_\k$ is \cite[Section~4,
Corollaire 1]{demazure;invariants-symetriques-entiers}.  The
multiplication rule follows from the product rule, Lemma
\ref{lemma;leibniz}, applied to $X=\pt$.
\end{proof}

The next result says that $H_T^*(X;\k)$ is in a natural way a
$\ca{D}_\k$-module.  This is analogous to a result of Brion in
equivariant Chow theory, \cite[Theorem
6.3]{brion;equivariant-chow-torus}.

\begin{theorem}\label{theorem;d-module}
The operations $\delta_\alpha$ and the action of $S_\k$ on
$H_T^*(X;\k)$ extend uniquely to an action of the algebra $\ca{D}_\k$.
This action is $H_G^*(X;\k)$-linear, it extends the action of the
twisted group algebra $S_\k[W]$, and it is contravariant with respect
to $G$-equivariant continuous maps and covariant with respect to
$G$-equivariant proper oriented maps.  In particular, the
characteristic homomorphism $S_\k\to H_T^*(X;\k)$ is
$\ca{D}_\k$-linear.
\end{theorem}

\begin{proof}
The basis elements $\partial_w$ of $\ca{D}_\k$ must act as in
\eqref{equation;partial}.  The multiplication rule given in
Proposition \ref{proposition;d} is then satisfied because of Lemma
\ref{lemma;leibniz}.  We conclude that the operations $\delta_\alpha$
generate a unique $\ca{D}_\k$-module structure on $H_T^*(X;\k)$.  It
follows from Lemma \ref{lemma;square} that
$\Delta(p_X^*(b)a)=p_X^*(b)\Delta(a)$ for all $\Delta\in\ca{D}_\k$,
$a\in H_T^*(X;\k)$ and $b\in H_G^*(X;\k)$, which shows that the
$\ca{D}_\k$-action is $H_G^*(X;\k)$-linear.  It follows from
Proposition \ref{proposition;divide-leibniz}\eqref{item;divide} that
the $\ca{D}_\k$-action extends the $S_\k[W]$-action.  Finally, the
functoriality follows from Lemma \ref{lemma;functorial}.
\end{proof}

By analogy with group rings, we define the \emph{augmentation left
ideal} of $\ca{D}_\k$ to be the annihilator of the constant polynomial
$1\in S_\k$,
$$I(\ca{D}_\k)=\{\,\Delta\in\ca{D}_\k\mid\Delta(1)=0\,\}.$$
Let $A$ be a left $\ca{D}_\k$-module.  We say an element of $A$ is
\emph{$\ca{D}_\k$-invariant} if it is annihilated by all operators in
$I(\ca{D}_\k)$ and we denote by $A^{I(\ca{D}_\k)}$ the set of
invariants.  The invariants are not a $\ca{D}_\k$-submodule of $A$,
but a submodule over the ring $(S^W)_\k$.

\begin{lemma}\label{lemma;annihilator}
\begin{enumerate}
\item\label{item;annihilator}
$I(\ca{D}_\k)$ is the left ideal of $\ca{D}_\k$ generated by the
operators $\partial_w$ for $w\ne1$.  It is a free left $S_\k$-module
with basis $(\partial_w)_{w\ne1}$.
\item\label{item;kill}
Let $A$ be a left $\ca{D}_\k$-module.  Then
$$
A^{I(\ca{D}_\k)}=\{\,a\in A\mid\text{$\partial_w(a)=0$ for all
$w\ne1$}\,\}\subset A^W.
$$
\end{enumerate}
\end{lemma}

\begin{proof}
\eqref{item;annihilator} follows from the fact that $\ca{D}_\k$ is a
free left $S_\k$-module with basis $(\partial_w)_{w\in W}$ and the
fact that $\partial_1=\id$ and $\partial_w(1)=0$ for $w\ne1$.  It
follows from \eqref{item;annihilator} that an element of $A$ is
$\ca{D}_\k$-invariant if and only if it is annihilated by all
$\partial_w$ for $w\ne1$.  By Proposition
\ref{proposition;divide-leibniz}\eqref{item;divide}, every
$\ca{D}_\k$-invariant element of $A$ is $W$-invariant.
\end{proof}

\begin{theorem}\label{theorem;annihilator}
We have inclusions 
$$
p_X^*\bigl(H_G^*(X;\k)\bigr)\subset H_T^*(X;\k)^{I(\ca{D}_\k)}\subset
H_T^*(X;\k)^W.
$$
\end{theorem}

\begin{proof}
The first inclusion follows from the $H_G^*(X;\k)$-linearity of the
$\ca{D}_\k$-action (Theorem \ref{theorem;d-module}) and the second
inclusion is a special case of Lemma
\ref{lemma;annihilator}\eqref{item;kill}.
\end{proof}

Without imposing conditions on the coefficient ring, this result
cannot be improved.  It is well-known that the map $p_X^*$ is not
always injective.  For instance, let $X$ be a point.  Then the kernel
of $p_X^*=p^*$ is the torsion submodule of $H^*(BG;\Z)$ (see e.g.\
\cite{feshbach;image-compact-lie} or Proposition
\ref{proposition;kernel-cokernel}\eqref{item;kernel} below), which is
nonzero for groups such as $G=\SO(n)$.  In Section
\ref{section;invariants} we will give some examples where the second
inclusion stated in the theorem is a strict inclusion.  Here is an
example where the first inclusion is strict.

\begin{example}\label{example;feshbach}
Let $X$ be a point, let $G=\Spin(11)$ or $\Spin(12)$, and let $\k=\Z$.
Then $H_T^*(X;\Z)=S$.  Put $S'=H_G^*(X;\Z)=H^*(BG;\Z)$.  If $p^*\colon
S'\to S^{I(\ca{D})}$ were surjective, we would have a ring isomorphism
$S'/\text{torsion}\cong S^{I(\ca{D})}$.  From the fact that $S$ is a
polynomial ring it follows easily that for every prime $l$ the ring
$S^{I(\ca{D})}\otimes_\Z\F_l$ has no nilpotents.  But Feshbach showed
in \cite{feshbach;image-compact-lie} that the ring
$(S'/\text{torsion})\otimes_\Z\F_2$ has nilpotents in degree $32$.
Thus $p^*\colon S'\to S^{I(\ca{D})}$ is not surjective.
\end{example}

In the next section, we shall give a criterion on the coefficient ring
$\k$ for $p_X^*$ to be an isomorphism from $H_G^*(X;\k)$ onto
$H_T^*(X;\k)^{I(\ca{D}_\k)}$.

\section{Equivariant cohomology and the torsion index}
\label{section;torsion}
  
As in the previous section, we fix a basis of the root system $R$ of
$(G,T)$.  We denote the set of positive roots by $R_+$, its
cardinality by $N$, and the longest element of $W$ by $w_0$.  As
before, $X$ denotes a topological $G$-space and $\k$ a commutative
ring.  The first result of this section says that the map in
cohomology induced by $p_X\colon X_T\to X_G$ ``almost'' has a left
inverse and ``almost'' surjects onto the Weyl invariants if
$H_T^*(X;\k)$ is torsion-free.  For a one-point space $X$, this result
can be found in \cite[Theorem~1.3]{totaro;torsion-index-spin}.

By a famous theorem of Borel,
\cite[Proposition~26.1]{borel;cohomologie-espaces-fibres}, the
characteristic homomorphism $c_G=i^*\colon S_\k\to H^*(G/T;\k)$ is
surjective over $\k=\Q$ and therefore it has finite cokernel over
$\k=\Z$.  Following \cite{grothendieck;torsion-homologique}, we define
the \emph{torsion index} $t(G)$ of $G$ to be the order of the cokernel
of $i^*\colon S^N\to H^{2N}(G/T;\Z)$.  It follows from
\eqref{equation;discriminant} that the torsion index divides the order
of the Weyl group.

\begin{proposition}\label{proposition;kernel-cokernel}
\begin{enumerate}
\item\label{item;kernel}
The kernel of $p_X^*\colon H_G^*(X;\k)\to H_T^*(X;\k)$ is annihilated
by $t(G)$.  If $t(G)$ is a unit in $\k$, then $p_X^*$ is split
injective.
\item\label{item;cokernel}
Assume that $X$ is compact, that $\k$ is a principal ideal domain of
characteristic $0$ and that the $\k$-module $H_T^*(X;\k)$ is
torsion-free.  Then the cokernel of
$$p_X^*\colon H_G^*(X;\k)\to H_T^*(X;\k)^W$$
is annihilated by $t(G)$.  If $t(G)$ is a unit in $\k$, $p_X^*$ is an
isomorphism onto $H_T^*(X;\k)^W$.
\end{enumerate}
\end{proposition}

\begin{proof}
Put $\phi=p_X^*$, $A=H_T^*(X;\k)$ and $B=H_G^*(X;\k)$.  By definition
of the torsion index, there exists a polynomial $u$ in $S_\k^N$ such
that $i^*(u)=t(G)\theta$, where $\theta\in H^{2N}(G/T;\Z)$ is the
orientation class.  Let $\bar{u}=\pr_T^*(u)\in A^{2N}$.  Then
$i_X^*(\bar{u})=i^*(u)=t(G)\theta$.  By \eqref{equation;integrate} and
by the naturality of $\partial_{w_0}$, this implies
\begin{multline*}
(p_Xi_X)^*p_{X,*}(\bar{u})=i_X^*\partial_{w_0,X}(\bar{u})
=\partial_{w_0,G}i_X^*(\bar{u})\\
=\partial_{w_0,G}i^*(u)=t(G)\partial_{w_0,G}(\theta)=t(G),
\end{multline*}
so $p_{X,*}(\bar{u})=t(G)$.  Define $\psi\colon A\to B$ by
$\psi(x)=p_{X,*}(\bar{u}x)$.  By the projection formula,
\begin{equation}\label{equation;torsion}
\psi(\phi(y))=p_{X,*}(\bar{u})y=t(G)y
\end{equation}
for all $y\in B$.  Hence $\phi(y)=0$ implies $t(G)y=0$, which shows
that the kernel of $\phi$ is killed by $t(G)$.  If $t(G)$ is
invertible in $\k$, then $t(G)^{-1}\psi$ is a left inverse of $\phi$,
so $\phi$ is split injective.  This proves \eqref{item;kernel}.  Now
make the assumptions stated in \eqref{item;cokernel}.  It follows from
\eqref{equation;torsion} that
$$\phi(\psi(\phi(y)))=\phi(t(G)y)=t(G)\phi(y)$$
for all $y\in B$, so
\begin{equation}\label{equation;image}
\phi(\psi(x))=t(G)x
\end{equation}
for all $x$ in the image of $\phi$.  Let $\K$ be the fraction field of
$\k$, $A_\K=A\otimes_\k\K$, $B_\K=B\otimes_\k\K$, and consider the
extended maps
$$
\phi_\K\colon B_\K\longto A_\K,\qquad\psi_\K\colon A_\K\longto B_\K.
$$
Because $X$ is compact, the universal coefficient theorem gives
$A_\K\cong H_T^*(X;\K)$ and $B_\K\cong H_G^*(X;\K)$.  Since $\K$ is a
field of characteristic $0$, it follows from \cite[Proposition
1(i)]{brion;equivariant-cohomology-intersection-theory} (or from
Theorem \ref{theorem;abelian} below) that $\phi_\K$ is an isomorphism
from $B_\K$ onto $(A_\K)^W$.  By \eqref{equation;image}, this implies
that $\phi_\K(\psi_\K(x))=t(G)x$ for all $x$ in $(A_\K)^W$.  But $A$
is by assumption torsion-free, so the natural map $A\to A_\K$ is
injective.  Hence the natural map $A^W\to(A_\K)^W$ is likewise
injective.  Therefore $\phi(\psi(x))=t(G)x$ for all $x$ in $A^W$,
which shows that the cokernel of $\phi$ is annihilated by $t(G)$.
Hence, if $t(G)$ is invertible, $\phi$ is surjective.
\end{proof}

Because of this proposition, we shall from now on work with a
coefficient ring in which the torsion index is invertible, such as
$\Z[t(G)^{-1}]$.  Our next result says that under this condition the
pullback diagram \eqref{equation;bundles} maps to a pushout diagram in
cohomology,
$$
\xymatrix{
H_T^*(X;\k)&H_G^*(X;\k)\ar[l]_{p_X^*}\\
H^*(BT;\k)\ar[u]^{\pr_T^*}&H^*(BG;\k)\ar[l]_{p^*}\ar[u]_{\pr_G^*}.
}
$$
The torsion index being invertible, there exists a polynomial
$\lie{S}$ in $S_\k^N$ such that
\begin{equation}\label{equation;fundamental}
i^*(\lie{S})=\theta\in H^{2N}(G/T;\k).
\end{equation}
This allows us to define the \emph{Schubert polynomials},
\begin{equation}\label{equation;schubert}
\lie{S}_w=\partial_{w^{-1}w_0}(\lie{S})\in S_\k^{l(w)},
\end{equation}
and their lifts to $T$-equivariant cohomology,
\begin{equation}\label{equation;classes}
\bar{\lie{S}}_w=\pr_T^*(\lie{S}_w)\in H_T^{2l(w)}(X;\k).
\end{equation}
Clearly $\lie{S}_{w_0}=\lie{S}$, where $w_0\in W$ is the longest
element.  It is shown in
\cite{demazure;invariants-symetriques-entiers} that $\lie{S}_1=1$, the
identity element of $S_\k^0$.  Hence $\bar{\lie{S}}_1=1\in
H_T^0(X;\k)$.  Moreover, by \cite[Corollaire
4]{demazure;invariants-symetriques-entiers}, the classes
$i^*(\lie{S}_w)$ form a basis (the Schubert basis) of the $\k$-module
$H^*(G/T;\k)$.  Thus we can define $\k$-linear maps
\begin{equation}\label{equation;split}
s\colon H^*(G/T;\k)\to S_\k,\qquad\bar{s}\colon H^*(G/T;\k)\to
H_T^*(X;\k)
\end{equation}
by setting $s(i^*(\lie{S}_w))=\lie{S}_w$ and
$\bar{s}(i^*(\lie{S}_w))=\bar{\lie{S}}_w$.

\begin{proposition}\label{proposition;restriction}
Assume that $t(G)$ is a unit in $\k$.  Choose a class $\lie{S}\in
S_\k^N$ as in \eqref{equation;fundamental}.  Define the classes
$\bar{\lie{S}}_w$ as in \eqref{equation;classes} and the map $\bar{s}$
as in \eqref{equation;split}.
\begin{enumerate}
\item\label{item;free}
The map $H^*(G/T;\k)\otimes_\k H_G^*(X;\k)\to H_T^*(X;\k)$ which sends
$(x,b)$ to $\bar{s}(x)p_X^*(b)$ is an isomorphism of
$H_G^*(X;\k)$-modules.  Hence $H_T^*(X;\k)$ is a free
$H_G^*(X;\k)$-module with basis $(\bar{\lie{S}}_w)_{w\in W}$.  In
particular, the rank of $H_T^*(X;\k)$ over $H_G^*(X;\k)$ is equal to
the cardinality of $W$ and the basis elements are of even
degree.
\item\label{item;tensor}
The map $S_\k\times H_G^*(X;\k)\to H_T^*(X;\k)$ which sends $(u,b)$ to
$\pr_T^*(u)p_X^*(b)$ induces an isomorphism of graded $\k$-algebras
$$
S_\k\otimes_{H^*(BG;\k)}H_G^*(X;\k)\longiso H_T^*(X;\k).
$$
\end{enumerate}
\end{proposition}

\begin{proof}
It follows from $i=\pr_T\circ i_X$ that
$i_X^*(\bar{\lie{S}}_w)=i^*(\lie{S}_w)$.  Thus $i_X^*$ is surjective
and $\bar{s}$ is a section of $i_X^*$.  Therefore the Leray-Hirsch
theorem applies to the fibre bundle $X_T\to X_G$, telling us that the
elements $\bar{\lie{S}}_w$ form a basis of the $H_G^*(X;\k)$-module
$H_T^*(X;\k)$.  This proves \eqref{item;free}.  The map $u\otimes
b\mapsto\pr_T^*(u)p_X^*(b)$ is multiplicative.  It sends the basis
element $\lie{S}_w\otimes1$ of $S_\k\otimes_{H^*(BG;\k)}H_G^*(X;\k)$
to the basis element $\bar{\lie{S}}_w$ of $H_T^*(X;\k)$ and is
therefore an algebra isomorphism.
\end{proof}

The next result is implicitly stated in
\cite{demazure;invariants-symetriques-entiers} and also follows from
\cite[Theorem~1.3]{totaro;torsion-index-spin}.  We will deduce it here
from Proposition \ref{proposition;restriction} by taking $X$ to be a
single point.  Let $M$ be a module over a commutative ring $\k$.  Let
$m\in M$ and let $I$ be an ideal of $\k$.  We call $m$ an
\emph{$I$-torsion element} if $I$ annihilates $m$.  We say $M$
\emph{has $I$-torsion} if $M$ has a nonzero $I$-torsion element.  If
$I$ is prime in $\k$ and $M$ has $I$-torsion, we call $I$ a
\emph{torsion prime} of $M$.

\begin{corollary}\label{corollary;bg-w}
Assume that $t(G)$ is a unit in $\k$.
\begin{enumerate}
\item\label{item;bg-w}
The maps $H^*(BG;\Z)\otimes_\Z\k\to H^*(BG;\k)$ and $p^*\colon
H^*(BG;\k)\to(S^W)_\k$ are isomorphisms.  Therefore $H^*(BG;\k)$ is a
polynomial algebra over $\k$ on generators of even degree.
\item\label{item;invariant}
The characteristic homomorphism $c_G=i^*\colon S_\k\to H^*(G/T;\k)$ is
surjective and its kernel is the ideal $(S^W_+)_\k$ of $S_\k$
generated by the Weyl-invariant elements of positive degree.  Hence
$c_G$ induces a $\ca{D}_\k$-linear isomorphism of $\k$-algebras
$$(S/S^W_+)_\k\longiso H^*(G/T;\k).$$
\end{enumerate}
\end{corollary}

\begin{proof}
By Proposition \ref{proposition;kernel-cokernel}\eqref{item;kernel},
the torsion submodule of $H^*(BG;\Z)$ is annihilated by the torsion
index $t(G)$.  Since $H^*(BG;\Z)$ is finitely generated in each
degree, this implies $\Tor^\Z(H^*(BG;\Z),\k)=0$, and therefore
$H^*(BG;\Z)\otimes_\Z\k\cong H^*(BG;\k)$ by the universal coefficient
theorem.  It follows that
$$
p^*(H^*(BG;\k))=p^*(H^*(BG;\Z))\otimes_\Z\k\subset(S^W)_\k.
$$
By Proposition \ref{proposition;restriction}\eqref{item;free}, applied
to a one-point space $X$, the algebra $S_\k=H^*(BT;\k)$ is a free
module with basis $(\lie{S}_w)_{w\in W}$ over the subalgebra
$p^*(H^*(BG;\k))$.  By \cite[Th\'eor\`eme
2(c)]{demazure;invariants-symetriques-entiers}, $S_\k$ is also a free
module with the same basis over the larger subalgebra $(S^W)_\k$.
This implies $p^*(H^*(BG;\k))$ is equal to $(S^W)_\k$.  By Proposition
\ref{proposition;kernel-cokernel}\eqref{item;kernel}, $p^*$ is
injective, and therefore $p^*\colon H^*(BG;\k)\to(S^W)_\k$ is an
isomorphism.  By \cite[Th\'eor\`eme
3]{demazure;invariants-symetriques-entiers}, it follows from this that
$H^*(BG;\k)$ is a polynomial algebra on even degree generators, which
proves \eqref{item;bg-w}.  By taking $X$ to be a space consisting of a
single point, we obtain from Proposition
\ref{proposition;restriction}\eqref{item;free}
$$
(S/S^W_+)_\k\cong S_\k\otimes_{(S^W)_\k}\k
\cong\bigl(H^*(G/T;\k)\otimes_\k(S^W)_\k\bigr)\otimes_{(S^W)_\k}\k\cong
H^*(G/T;\k).
$$
This isomorphism is $\ca{D}_\k$-linear, because $i^*$ is.
\end{proof}

We shall from now on identify $H^*(BG;\k)$ with $(S^W)_\k$ via the
isomorphism $p^*$ (when the torsion index is invertible in $\k$).

\begin{corollary}\label{corollary;d-module}
Assume that $t(G)$ is a unit in $\k$.  The isomorphism
\begin{equation}\label{equation;tensor}
S_\k\otimes_{(S^W)_\k}H_G^*(X;\k)\cong H_T^*(X;\k)
\end{equation}
of Proposition {\rm\ref{proposition;restriction}}\eqref{item;tensor}
is $\ca{D}_\k$-linear, where we let the operators $\Delta\in\ca{D}_\k$
act on the left-hand side by $\Delta(u\otimes b)=\Delta(u)\otimes b$.
\end{corollary}

\begin{proof}
By Theorem \ref{theorem;d-module}, the action of $\ca{D}_\k$ on $S_\k$
is linear over $H^*(BG;\k)=(S^W)_\k$, which implies that the operation
$\Delta(u\otimes b)=\Delta(u)\otimes b$ is well-defined.  Again by
Theorem \ref{theorem;d-module}, the action of $\ca{D}_\k$ on
$H_T^*(X;\k)$ is $H_G^*(X;\k)$-linear, which implies that the
isomorphism $f(u\otimes b)=\pr_T^*(u)p_X^*(b)$ satisfies
$f\Delta(u\otimes b))=f(\Delta(u)\otimes b)=\Delta f(u\otimes b)$.
\end{proof}

We recall from \cite{demazure;invariants-symetriques-entiers} that the
prime factors of the torsion index are precisely the torsion primes of
the fundamental group of $G$ plus a short list of primes for each of
the simple ideals of the Lie algebra of $G$.  This list is as follows:
$2$ for types $\B_n$ ($n\ge3$), $\D_n$ ($n\ge4$) and $\G_2$; $2$ and
$3$ for types $\E_6$, $\E_7$ and $\F_4$; and $2$, $3$ and $5$ for type
$\E_8$.  (The simply connected groups of type $\A_n$ and $\C_n$ have
torsion index $1$ for all $n\ge1$.)  The proof of this fact is largely
given in \cite{borel;sous-groupes-commutatifs} and relies on the
observation (also used in the proof of Corollary \ref{corollary;bg-w})
that the prime factors of the torsion index are precisely the torsion
primes of $H^*(BG;\Z)$, which in turn are the same as the torsion
primes of $H^*(G;\Z)$.  Thus $t(G)$ is invertible in $\k$ if and only
if the torsion primes of $H^*(G;\Z)$ are invertible in $\k$.

Proposition \ref{proposition;restriction}\eqref{item;tensor} shows
that the $T$-equivariant cohomology of a $G$-space is determined by
its $G$-equivariant cohomology for any coefficient ring in which the
torsion index is invertible.  Conversely, Theorem \ref{theorem;d}
below shows that it is possible to express $G$-equivariant cohomology
in terms of $T$-equivariant cohomology.  This result is a consequence
of the following theorem of Demazure.

\begin{theorem}[{\cite[Th\'eor\`eme
2]{demazure;invariants-symetriques-entiers}}]
\label{theorem;demazure}
Assume that $t(G)$ is a unit in $\k$.
\begin{enumerate}
\item\label{item;polynomial-basis}
Choose a family of Schubert polynomials $(\lie{S}_w)_{w\in W}$ as in
\eqref{equation;schubert}.  These polynomials form a basis of the
$(S^W)_\k$-module $S_\k$.
\item\label{item;operator-basis}
$\ca{D}_\k=\End_{(S^W)_\k}(S_\k)$.
\end{enumerate}
\end{theorem}

It follows from this theorem that $\ca{D}_\k$ is isomorphic to the
algebra of $\abs{W}\times\abs{W}$-matrices over the ring $(S^W)_\k$.
In particular, the centre of $\ca{D}_\k$ is equal to $(S^W)_\k$.

\begin{theorem}\label{theorem;d}
Assume that $t(G)$ is a unit in $\k$.
\begin{enumerate}
\item\label{item;d}
The map $p_X^*$ is an isomorphism from $H_G^*(X;\k)$ onto
$H_T^*(X;\k)^{I(\ca{D}_\k)}$.
\item\label{item;cohomology-operator-basis}
The family of operators $(\partial_w)_{w\in W}$ is a basis of the
$H_T^*(X;\k)$-module of all $H_G^*(X;\k)$-linear endomorphisms of
$H_T^*(X;\k)$.
\end{enumerate}
\end{theorem}

\begin{proof}
Let
$\lie{A}=\ca{D}_\k$-$\lie{Mod}$ and $\lie{B}=(S_\k)^W$-$\lie{Mod}$ be
the categories of (left) modules over the rings $\ca{D}_\k$, resp.\
$(S_\k)^W$.  By Theorem
\ref{theorem;demazure}\eqref{item;polynomial-basis}, the module $S_\k$
is a progenerator of the category $\lie{B}$.  Hence, by the first
Morita equivalence theorem (see e.g.\
\cite[\S~18]{lam;modules-rings}), the functor
$\lie{G}\colon\lie{B}\to\lie{A}$ defined by
$$B\longmapsto S_\k\otimes_{(S^W)_\k}B$$
is an equivalence with inverse $\lie{F}\colon\lie{A}\to\lie{B}$ given
by
$$A\longmapsto\Hom_{\ca{D}_\k}(S_\k,A).$$
We assert that $\lie{F}$ is naturally isomorphic to the functor
$\lie{I}\colon\lie{A}\to\lie{B}$ given
by
$$A\longmapsto A^{I(\ca{D}_\k)}.$$
Indeed, consider the natural $(S_\k)^W$-linear map
$\Phi_A\colon\Hom_{\ca{D}_\k}(S_\k,A)\to A$ defined by
$\Phi_A(f)=f(1)$.  The map $\Phi_A$ is injective, because
$\ca{D}_\k\supset S_\k$.  Its image is $A^{I(\ca{D}_\k)}$, because
$I(\ca{D}_\k)$ is the annihilator of $1\in S_\k$.  We conclude that
$\Phi$ is a natural isomorphism from $\lie{F}$ to $\lie{I}$.  Now
consider the $\ca{D}_\k$-module $A=H_T^*(X;\k)$ and the
$(S^W)_\k$-module $B=H_G^*(X;\k)$.  By Proposition
\ref{proposition;restriction}\eqref{item;tensor}, $A\cong\lie{G}(B)$.
Hence $B\cong\lie{F}(A)\cong\lie{I}(A)=A^{I(\ca{D}_\k)}$, which proves
\eqref{item;d}.  The isomorphism \eqref{equation;tensor} gives an
isomorphism of $H_T^*(X;\k)$-modules
$$
\End_{H_G^*(X;\k)}(H_T^*(X;\k))
\cong\End_{(S^W)_\k}(S_\k)\otimes_{(S^W)_\k}H_G^*(X;\k).
$$
It now follows from Proposition \ref{proposition;d} and Theorem
\ref{theorem;demazure}\eqref{item;operator-basis} that the
$\partial_w$ form a basis of the $H_T^*(X;\k)$-module
$\End_{H_G^*(X;\k)}(H_T^*(X;\k))$.
\end{proof}

\begin{example}\label{example;d}
Let $X$ be a point.  It follows from Corollary \ref{corollary;bg-w}
and Theorem \ref{theorem;d} that for any $\k$ in which $t(G)$ is
invertible
$$
(S^W)_\k\cong H^*(BG;\k)\cong(S_\k)^{I(\ca{D}_\k)}.
$$
We shall see in Corollary \ref{corollary;abelian} that
$(S_\k)^{I(\ca{D}_\k)}=(S_\k)^W$ under a mild condition on the root
data.
\end{example}

In general, $H_G^*(X;\k)$ is a proper submodule of $H_T^*(X;\k)^W$.
(See Section \ref{section;invariants} for examples.)  Nevertheless,
there are two useful sufficient criteria for the inclusion to be an
equality, stated in Theorem \ref{theorem;abelian} below.  The second
criterion (that the cardinality of the Weyl group be invertible in
$\k$) is well-known, but it seems to us of some interest to deduce it
from Theorem \ref{theorem;d}.

\begin{lemma}\label{lemma;left-inverse}
The notation and the hypotheses are as in Theorem
{\rm\ref{theorem;demazure}}.  Define $\psi\in\ca{D}_\k$ by
$\psi(u)=\partial_{w_0}(\lie{S}_{w_0}u)$ for $u\in S_\k$.  Let $J$ be
the left ideal of $\ca{D}_\k$ generated by $\psi$ and by all elements
$1-w$ for $w\in W$.  Let $A$ be a left $\ca{D}_\k$-module.
\begin{enumerate}
\item\label{item;left-inverse}
$\psi$ projects $A$ onto the $(S^W)_\k$-submodule $A^{I(\ca{D}_\k)}$.
\item\label{item;left-ideal}
$A^W=A^{I(\ca{D}_\k)}\oplus A^J$.  Hence $A^W=A^{I(\ca{D}_\k)}$ if and
only if $A^J=0$.
\end{enumerate}
\end{lemma}

\begin{proof}
It follows from \eqref{equation;long-symmetric} that $\partial_{w_0}$
maps $S_\k$ to $(S^W)_\k$ and hence that $\partial_w\partial_{w_0}=0$
for all $w\ne1$.  Therefore $\partial_w(\psi(a))=0$ for all $a\in A$
and all $w\ne1$, i.e.\ $\psi$ maps $A$ to $A^{I(\ca{D}_\k)}$.  If
$a\in A^{I(\ca{D}_\k)}$, then by Lemma \ref{lemma;leibniz}
$$
\psi(a)=\partial_{w_0}(\lie{S}_{w_0}a)
=\partial_{w_0}(\lie{S}_{w_0})a=a,
$$
since $\partial_{w_0}(\lie{S}_{w_0})=\lie{S}_1=1$.  This proves that
$\psi(A)=A^{I(\ca{D}_\k)}$ and $\psi^2=\id_A$, which establishes
\eqref{item;left-inverse}.  It follows from \eqref{item;left-inverse}
that $A$ is the direct sum of the $(S^W)_\k$-submodules
$A^{I(\ca{D}_\k)}$ and $\ker(\psi)$.  Moreover, it follows from
\eqref{item;left-inverse} and from Lemma
\ref{lemma;annihilator}\eqref{item;kill} that $\psi$ maps $A^W$ into
itself.  Therefore $A^W$ is the direct sum of $A^{I(\ca{D}_\k)}$ and
$A^W\cap\ker(\psi)=A^J$.
\end{proof}

For certain groups $G$ it turns out that there exists an ideal
$J(S_\k)$ of $S_\k$ such that, for all $\ca{D}_\k$-modules $A$, we
have $A^W=A^{I(\ca{D}_\k)}$ if and only if $A^{J(S_\k)}=0$.  (See
Section \ref{section;invariants} for examples.)  We do not know if
this is true for all groups.

Recall that $d=\prod_{\alpha\in R_+}\alpha\in S^N$ denotes the
discriminant of $G$.  We shall denote its image in $S_\k^N$ also by
$d$.

\begin{lemma}\label{lemma;discriminant}
The notation and the hypotheses are as in Lemma
{\rm\ref{lemma;left-inverse}}.
\begin{enumerate}
\item\label{item;discriminant}
The discriminant satisfies the identity
$$
d(1-\psi)=\sum_{w\in W}\det(w)w(\lie{S}_{w_0})(1-w).
$$
In particular, $d\in J$.  Hence $A^W=A^{I(\ca{D}_\k)}$ if $d$ is not a
zero divisor in $A$.
\item\label{item;weyl-invertible}
Assume that $\abs{W}$ is a unit in $\k$.  Choose the top-degree
Schubert polynomial to be $\lie{S}_{w_0}=\abs{W}^{-1}d$.  Then $\psi$
is the projection map onto the $W$-invariants,
$$
\psi=\frac1{\abs{W}}\sum_{w\in W}w.
$$
Hence $A^W=A^{I(\ca{D}_\k)}$.
\end{enumerate}
\end{lemma}

\begin{proof}
It follows from \eqref{equation;long-symmetric} that for all $u\in
S_\k$
\begin{multline*}
d\psi(u)=d\partial_{w_0}(\lie{S}_{w_0}u)=\sum_{w\in
W}\det(w)w(\lie{S}_{w_0}u)\\
=\sum_{w\in W}\det(w)w(\lie{S}_{w_0})u+\sum_{w\in
W}\det(w)w(\lie{S}_{w_0})(w(u)-u)\\
=d\partial_{w_0}(\lie{S}_{w_0})u+\sum_{w\in
W}\det(w)w(\lie{S}_{w_0})(w-1)(u)\\
=du+\sum_{w\in W}\det(w)w(\lie{S}_{w_0})(w-1)(u).
\end{multline*}
This identity shows that $d\in J$.  Hence, if $d$ is not a zero
divisor in $A$, then $A^J=0$, so $A^W=A^{I(\ca{D}_\k)}$ by Lemma
\ref{lemma;left-inverse}\eqref{item;left-ideal}.  This proves
\eqref{item;discriminant}.  Assume $\abs{W}$ is a unit in $\k$.
Comparing \eqref{equation;discriminant} with
\eqref{equation;fundamental} we see that $\lie{S}_{w_0}=\abs{W}^{-1}d$
is a valid choice for the top Schubert polynomial.  It follows from
\eqref{equation;antisymmetry} and the fact that $d$ is antisymmetric
that
$$
\partial_{w_0}(du)=\frac1d\sum_{w\in W}\det(w)w(du)=\sum_{w\in W}w(u)
$$
for all $u\in S_\k$.  Hence
$\psi(u)=\partial_{w_0}(\lie{S}_{w_0}u)=\abs{W}^{-1}\sum_{w\in
W}w(u)$.  Lemma \ref{lemma;left-inverse}\eqref{item;left-inverse} now
implies that $A^W=A^{I(\ca{D}_\k)}$.
\end{proof}

\begin{theorem}\label{theorem;abelian}
Assume that $t(G)$ is a unit in $\k$.  Assume either that $d$ is not a
zero divisor in the $S_\k$-module $H_T^*(X;\k)$ or that $\abs{W}$ is a
unit in $\k$.  Then $p_X^*$ is an isomorphism from $H_G^*(X;\k)$ onto
$H_T^*(X;\k)^W$.
\end{theorem}

\begin{proof}
This follows immediately from Theorem \ref{theorem;d} and Lemma
\ref{lemma;discriminant}.
\end{proof}

\begin{corollary}\label{corollary;abelian-torsion}
Assume that $t(G)$ is a unit in $\k$ and that the restriction map
$H_T^*(X;\k)\to H_T^*(X^T;\k)$ is injective.  In addition, if $G$ has
a root $\alpha$ such that $\alpha/2\in\X(T)$ {\rm(}i.e.\ if $G$
contains a direct factor isomorphic to $\U(n,\H)$ for some
$n\ge1${\rm)}, assume that $H^*(X^T;\k)$ has no $2$-torsion.  Then
$p_X^*\colon H_G^*(X;\k)\to H_T^*(X;\k)^W$ is an isomorphism.
\end{corollary}

\begin{proof}
Let $\alpha$ be a root of $G$.  If $\alpha/2$ is not a character of
$T$, put $\xi_1=\alpha$.  If $\alpha/2$ is a character, put
$\xi_1=\alpha/2$.  In either case, $\Z\xi_1$ is a direct summand of
the character group $\X(T)$.  Select characters $\xi_2$,
$\xi_3$,\dots, $\xi_r$ such that $(\xi_1,\xi_2,\dots,\xi_r)$ is a
basis of $\X(T)$.  Let $\eta_j=\xi_j\otimes1_\k$ be the image of
$\xi_j$ in $\X(T)_\k=\X(T)\otimes_\Z\k$.  Then
$(\eta_1,\eta_2,\dots,\eta_r)$ is a basis of $\X(T)_\k$, so
$H_T^*(X^T;\k)\cong H^*(X^T;\k)\otimes_\k S_\k$ is a polynomial
algebra over $H^*(X^T;\k)$ in the variables $\eta_1$, $\eta_2$,\dots,
$\eta_r$.  The element $\alpha\otimes1_\k$ is equal to $\eta_1$ if
$\alpha/2\not\in\X(T)$ and to $2\eta_1$ if $\alpha/2\in\X(T)$.  In
either case, $\alpha\otimes1_\k$ is not a zero divisor in
$H_T^*(X^T;\k)$.  This implies that $d$ is not a zero divisor in the
submodule $H_T^*(X;\k)$ of $H_T^*(X^T;\k)$.  The corollary now follows
from Theorem \ref{theorem;abelian}.
\end{proof}

The next corollary is a result stated in
\cite{demazure;invariants-symetriques-entiers} under the marginally
stronger assumption that $2$ be invertible in $\k$.  (Note, however,
that Demazure's proof works equally well under our weaker
assumptions.)

\begin{corollary}[{\cite[Corollaire au Th\'eor\`eme
2]{demazure;invariants-symetriques-entiers}}]
\label{corollary;abelian}
Assume that $t(G)$ is a unit in $\k$.  In addition, if $G$ has a root
$\alpha$ such that $\alpha/2\in\X(T)$, assume that $2$ is not a zero
divisor in $\k$.  Then $p^*\colon H^*(BG;\k)\to H^*(BT;\k)^W=(S_\k)^W$
is an isomorphism.  Thus
$$
(S_\k)^W\cong H^*(BG;\k)\cong H^*(BG;\Z)\otimes_\Z\k\cong(S^W)_\k.
$$
\end{corollary}

\begin{proof}
The first statement follows from Corollary
\ref{corollary;abelian-torsion}, applied to $X=\rm{pt}$.  The other
isomorphisms were established in Corollary \ref{corollary;bg-w}.
\end{proof}

\section{Orthogonality}\label{section;orthogonal}

In this section we obtain a duality result based on material in
\cite[Chapter V]{macdonald;notes-schubert-polynomials}.  As in the
previous sections, we let $X$ be a topological $G$-space and $\k$ a
commutative ring.  We fix a basis of the root system of $(G,T)$ and
relative to this basis define the operators $\partial_w$ in Demazure's
algebra $\ca{D}_\k$ as in \eqref{equation;partial}.  Define
$$\ca{B}_X(a_1,a_2)=\partial_{w_0}(a_1a_2)$$
for $a_1$, $a_2\in H_T^*(X;\k)$.  It follows from the relations
\eqref{equation;relation} that $\Delta\partial_{w_0}=0$ for all
$\Delta$ in the augmentation left ideal $I(\ca{D}_\k)$, so this rule
defines a pairing
$$
\ca{B}_X\colon H_T^*(X;\k)\times H_T^*(X;\k)\to
H_T^{*-2N}(X;\k)^{I(\ca{D}_\k)}.
$$
This pairing is graded symmetric in the sense that
$$\ca{B}_X(a_1,a_2)=(-1)^{k_1k_2}\ca{B}_X(a_2,a_1)$$
if $a_1$ is of degree $k_1$ and $a_2$ is of degree $k_2$.  It follows
from the product rule, Lemma \ref{lemma;leibniz}, that
$H_T^*(X;\k)^{I(\ca{D}_\k)}$ is a subalgebra of $H_T^*(X;\k)$ and that
the pairing is bilinear over $H_T^*(X;\k)^{I(\ca{D}_\k)}$ in the sense
that it is additive in each variable and
$$
\ca{B}_X(ba_1,a_2)=b\ca{B}_X(a_1,a_2),\qquad
\ca{B}_X(a_1,a_2b)=\ca{B}_X(a_1,a_2)b
$$
for all $b\in H_T^*(X;\k)^{I(\ca{D}_\k)}$.  It follows from the
naturality of $\partial_{w_0}$ that
\begin{equation}\label{equation;natural-pairing}
\ca{B}_X(f_T^*(a_1),f_T^*(a_2))=f_T^*\ca{B}_Y(a_1,a_2)
\end{equation} 
for $a_1$, $a_2\in H_T^*(Y;\k)$, where $f_T\colon X_T\to Y_T$ is the
map induced by any $G$-equivariant continuous map $f\colon X\to Y$.
The pairing may be degenerate, but nevertheless the operators
$\partial_w$ and $w$ have adjoints.

\begin{lemma}\label{lemma;adjoint}
Let $w\in W$ and $a_1$, $a_2\in H_T^*(X;\k)$.  Then
\begin{enumerate}
\item\label{item;adjoint-partial}
$\ca{B}_X(\partial_w(a_1),a_2)=\ca{B}_X(a_1,\partial_{w^{-1}}(a_2))$,
\item\label{item;adjoint-w}
$\ca{B}_X(w(a_1),a_2)=\det(w)\ca{B}_X(a_1,w^{-1}(a_2))$.
\end{enumerate}
\end{lemma}

\begin{proof}
To prove \eqref{item;adjoint-partial}, it is enough to show that
$\ca{B}_X(\partial_{s_\alpha}a_1,a_2)=\ca{B}_X(a_1,\partial_{s_\alpha}a_2)$
for all simple roots $\alpha$.  By Proposition
\ref{proposition;divide-leibniz}\eqref{item;leibniz}, Lemma
\ref{lemma;reflect} and \eqref{equation;relation},
$$
\ca{B}_X(\partial_{s_\alpha}(a_1),a_2)
=\partial_{w_0s_\alpha}\partial_{s_\alpha}(\partial_{s_\alpha}(a_1)a_2)
=\partial_{w_0s_\alpha}(\partial_{s_\alpha}(a_1)\partial_{s_\alpha}(a_2)).
$$
This expression is graded symmetric in $a_1$ and $a_2$, so assuming
that $a_1$ is of degree $k_1$ and $a_2$ is of degree $k_2$ we have
$$
\ca{B}_X(\partial_{s_\alpha}(a_1),a_2)
=(-1)^{k_1k_2}\ca{B}_X(\partial_{s_\alpha}(a_2),a_1)
=\ca{B}_X(a_1,\partial_{s_\alpha}(a_2)),
$$
which proves \eqref{item;adjoint-partial}.  From this
\eqref{item;adjoint-w} follows by using Proposition
\ref{proposition;divide-leibniz}\eqref{item;divide}.
\end{proof}

When $G$ is the unitary group, the lifted Schubert classes
$\bar{\lie{S}}_w$ have a dual basis.  This leads to a relative form of
Poincar\'e duality for the fibre bundle $X_T\to X_G$.

\begin{proposition}\label{proposition;orthogonal}
Let $G=\U(l+1)$.  Then the pairing
$$
\inner{\mult,\mult}\colon H_T^*(X;\k)\times H_T^*(X;\k)\to
H_G^{*-2N}(X;\k)
$$
defined by $\inner{a_1,a_2}=p_{X,*}(a_1a_2)$ is nonsingular.  Hence
$$
H_T^*(X;\k)\cong
\Hom_{H_G^*(X;\k)}\bigl(H_T^*(X;\k),H_G^{*-2N}(X;\k)\bigr)
$$
as graded left $H_G^*(X;\k)$-modules.
\end{proposition}

\begin{proof}
Let $\lie{S}=\eps_1^l\eps_2^{l-1}\cdots\eps_l$ be the top Schubert
polynomial, where $\eps_k$ is the character of the diagonal maximal
torus of $G$ defined by
\begin{equation}\label{equation;character}
\eps_k\begin{pmatrix}t_1&&\\&\zerodots&\\&&t_{l+1}\end{pmatrix}=t_k.
\end{equation}
For each $w\in W$, define $\bar{\lie{S}}_w\in H_T^{2l(w)}(X;\k)$ as in
\eqref{equation;classes} and define
$$\bar{\lie{S}}^w=\det(ww_0)w_0(\bar{\lie{S}}_{ww_0}).$$
Applying \eqref{equation;natural-pairing} to the constant map
$X\to\pt$, we find
$$
\ca{B}_X(\bar{\lie{S}}_w,\bar{\lie{S}}^{w'})
=\ca{B}_X(\pr_T^*(\lie{S}_w),\pr_T^*(\lie{S}^{w'}))
=\pr_T^*\ca{B}_\pt(\lie{S}_w,\lie{S}^{w'})=\delta_{w,w'}
$$
for all $w$, $w'\in W$, where the last equality follows from
\cite[V(5.5)]{macdonald;notes-schubert-polynomials}.  This shows that
the pairing $\ca{B}_X$ is nonsingular.  Since
$p_X^*(\inner{a_1,a_2})=\ca{B}_X(a_1,a_2)$ and $p_X^*$ is an
isomorphism from $H_G^*(X;\k))$ onto $H_T^*(X;\k)^{I(\ca{D}_\k)}$
(Theorem \ref{theorem;d}), the pairing $\inner{\mult,\mult}$ is
nonsingular as well.
\end{proof}

This result is presumably true for arbitrary $G$, provided that the
torsion index is a unit in $\k$, but we do not know if the Schubert
polynomials can be chosen so that
$\ca{B}_\pt(\lie{S}_w,\lie{S}^{w'})=\delta_{w,w'}$.

\section{Weyl invariants for groups of low rank}
\label{section;invariants}

In this section we calculate, for a few groups $G$ of low rank,
$H_T^*(X;\k)^W$ as a module over $H_G^*(X;\k)$.  We do this by using
Proposition \ref{proposition;restriction}, which says that
$H_T^*(X;\k)$ is a free $H_G^*(X;\k)$-module with basis
$(\bar{\lie{S}}_w)_{w\in W}$.  (We consider only coefficient rings
$\k$ in which the torsion index of $G$ is invertible.)  Finding the
Weyl invariants is then a matter of calculating the matrices of the
simple reflections with respect to this basis and solving linear
equations.  The results are summarized in the table below.

Here we write $A=H_T^*(X;\k)$ and $B=H_G^*(X;\k)$ and identify $B$
with the subalgebra $p_X^*(B)$ of $A$.  We write $B^I$ for the
submodule of $B$ annihilated by an ideal $I$ of $(S^W)_\k$.  To save
space, instead of $A^W$ we give the $B$-submodule $A^J$, where $J$ is
the ideal of $\ca{D}_\k$ defined in Lemma \ref{lemma;left-inverse}.
(Recall that $A^W=B\oplus A^J$ by Lemma \ref{lemma;left-inverse} and
Theorem \ref{theorem;d}.)

The rightmost column of the table shows a certain ideal $J(S_\k)$ of
$S_\k$, which has the property that $A^W=B$ if and only if
$A^{J(S_\k)}=0$.  This ideal is found as follows.  Upon examining the
third column one sees that, for each of the groups we consider, there
exists an ideal $K$ of $(S^W)_\k$ such that $A^W=B$ if and only if
$B^K=0$.  For instance, $K=(2,p_1)$ for $G=\U(2)$, $K=(2,p_1p_2+p_3)$
for $G=\U(3)$, and $K=(2)$ for $G=\U(2,\H)$.  We define $J(S_\k)$ to
be the ideal of $S_\k$ generated by $K$.  This has the desired
property, because $S_\k$, being a free $(S^W)_\k$-module, is
faithfully flat over $(S^W)_\k$.  It follows from this that
$K=(S^W)_\k\cap J(S_\k)$ and $A^{J(S_\k)}=S_\k\otimes_{(S^W)_\k}B^K$.
Thus $B^K=0$ if and only if $A^{J(S_\k)}=0$.  Observe that $d\in
J(S_\k)$ in all cases.
$$
\begin{tabular}{@{}lccc@{}}
\toprule
$G$&$t(G)$&$A^J$&$J(S_\k)$\\
\midrule
$\U(2)$&$1$&$B^{(2,p_1)}\mult\bar{\lie{S}}_{w_0}$&$(2,d)$\\
$\U(3)$&$1$&$B^{(2,p_1p_2+p_3)}
\mult(\bar{\lie{S}}_{w_0}+p_1\bar{\lie{S}}_{s_2s_1}
+p_2\bar{\lie{S}}_{s_2}+p_1^2\bar{\lie{S}}_{s_1})$&$(2,d)$\\
$\SU(2)$&$1$&$B^{(2)}\mult\bar{\lie{S}}_{w_0}$&$(2)$\\
$\SU(3)$&$1$&$B^{(2,q_3)}
\mult(\bar{\lie{S}}_{w_0}+q_2\bar{\lie{S}}_{s_2})$&$(2,d)$\\
$\SO(3)$&$2$&$0$&$(1)$\\
$\PSU(3)$&$3$&$B^{(2,q_3)}
\mult(\bar{\lie{S}}_{w_0}+q_2\bar{\lie{S}}_{s_2})$&$(2,d)$\\
$\U(2,\H)$&$1$&$B^{(2,p_1)}\mult\bar{\lie{S}}_{w_0}\oplus
B^{(2)}\mult\bar{\lie{S}}_{s_2s_1s_2}\oplus
B^{(2)}\mult\bar{\lie{S}}_{s_1s_2}\qquad$&$(2)$\\
&&$\qquad\qquad\qquad\qquad\oplus
B^{(2,p_1)}\mult\bar{\lie{S}}_{s_2s_1}\oplus
B^{(2)}\mult\bar{\lie{S}}_{s_2}$\\
\bottomrule
\end{tabular}
$$

We use the notation of \cite[Chapitres IV--VI, Planches
I--IX]{bourbaki;groupes-algebres} for roots and weights.  The other
notations, as well as a few computational details, are explained in
the discussion below.

\subsection*{$\U(2)$ and $\U(3)$}

Let $G=\U(l+1)$, let $T$ be the diagonal maximal torus and let
$\eps_1$, $\eps_2$,\dots, $\eps_{l+1}$ be the basis elements of
$\X(T)$ defined in \eqref{equation;character}.  Let
$\lie{S}=\eps_1^l\eps_2^{l-1}\cdots\eps_l$.  By
\cite[Section~5]{demazure;invariants-symetriques-entiers},
$i^*(\lie{S})=i^*(d)/(l+1)!\in H^{2N}(G/T;\Z)$, so $t(G)=1$ and the
coefficient ring $\k$ is arbitrary.  The Schubert polynomials
$\lie{S}_w=\partial_{w^{-1}w_0}(\lie{S})$ form a basis of the
$S^W$-module $S$ and the algebra $S^W$ is freely generated by the
elementary symmetric polynomials $p_1$, $p_2$,\dots, $p_{l+1}$.  For
$\U(2)$ the Schubert polynomials are $\lie{S}_s=\lie{S}=\eps_1$ and
$\lie{S}_1=1$, where $s$ denotes the generator of $W$, which acts on
$\X(T)$ by interchanging $\eps_1$ and $\eps_2$.  The matrix of $s$
with respect to the basis $(\lie{S}_s,1)$ is
$$
M=\begin{pmatrix}-1&0\\p_1&1\end{pmatrix}.
$$
This is also the matrix of $s$ with respect to the basis
$(\bar{\lie{S}}_s,1)$ of the $B$-module $A$.  Thus an element of $A$
of the form $c_1\bar{\lie{S}}_s+c_2$ with $c_1$, $c_2\in B$ is
$W$-invariant if and only if
$M\bigl(\!\begin{smallmatrix}c_1\\c_2\end{smallmatrix}\!\bigr)
=\bigl(\!\begin{smallmatrix}c_1\\c_2\end{smallmatrix}\!\bigr)$, which
is the case if and only if $2c_1=p_1c_1=0$.  The upshot is
$A^W=B^{(2,p_1)}\mult\bar{\lie{S}}_s\oplus B$.  For $\U(3)$ the
Schubert polynomials are
\begin{gather*}
\lie{S}_{w_0}=\lie{S}=\eps_1^2\eps_2,\quad
\lie{S}_{s_1s_2}=\eps_1\eps_2,\quad\lie{S}_{s_2s_1}=\eps_1^2,\\
\lie{S}_{s_2}=\eps_1+\eps_2,\quad
\lie{S}_{s_1}=\eps_1,\quad\lie{S}_1=1,
\end{gather*}
where $s_1$ and $s_2$ denote the reflections in the roots
$\alpha_1=\eps_1-\eps_2$ and $\alpha_2=\eps_2-\eps_3$.  The matrices
of $s_1$ and $s_2$ with respect to the Schubert basis
$(\lie{S}_w)_{w\in W}$ are
$$
M_1=
\begin{pmatrix}
-1&0&0&0&0&0\\
p_1&1&-1&0&0&0\\
0&0&-1&0&0&0\\
0&0&p_1&1&1&0\\
0&0&0&0&-1&0\\
-p_3&0&-p_2&0&0&1
\end{pmatrix},
\qquad
M_2=
\begin{pmatrix}
-1&0&0&0&0&0\\
0&-1&0&0&0&0\\
0&-1&1&0&0&0\\
0&0&0&-1&0&0\\
p_2&p_1&0&1&1&0\\
-p_3&0&0&p_1&0&1
\end{pmatrix}.
$$
To compute $A^W$ one now solves the linear equations
$M_1\vec{c}=M_2\vec{c}=\vec{c}$ for a column vector
$\vec{c}={}^t(c_1,c_2,c_3,c_4,c_5,c_6)$ with entries in $B$.

\subsection*{$\SU(2)$ and $\SU(3)$}

Similar results for $G=\SU(l+1)$ are obtained by restricting the
computations done for $\U(l+1)$ to the maximal torus of $G$, i.e.\ by
replacing $\eps_k$ with $\eps_k-\eps_0$, where
$\eps_0=(\eps_1+\eps_2+\cdots+\eps_{l+1})/(l+1)$.  We denote by $q_k$
the restriction to $\ca{X}(T)$ of the elementary symmetric polynomial
$p_k$.  (In particular, $q_1=0$.)  By
\cite[Section~5]{demazure;invariants-symetriques-entiers}, we have
$t(G)=1$, so again we can use any coefficient ring $\k$.

\subsection*{$\SO(3)$ and $\PSU(3)$}

Let $G=\PSU(l+1)$, the adjoint form of $\tilde{G}=\SU(l+1)$, and let
$T$ be the image in $G$ of the diagonal maximal torus $\tilde{T}$ of
$\tilde{G}$.  The character group $\X(T)$ has index $l+1$ in
$\X(\tilde{T})$ and, by
\cite[Section~2]{totaro;torsion-index-8-other}, $t(G)=l+1$.
Accordingly, we can let $\k$ be any commutative ring in which $l+1$ is
invertible.  Then the canonical map $S(\X(T))_\k\to
S(\X(\tilde{T}))_\k$ is an isomorphism, so the calculations done above
for $\tilde{G}$ carry over directly to $G$.

\subsection*{$\U(2,\H)$}

Let $G=\U(2,\H)$.  We regard $\U(2)$ as a subgroup of $G$; then the
diagonal maximal torus $T$ of $\U(2)$ is also maximal in $G$.  By
\cite[Section~5]{demazure;invariants-symetriques-entiers}, $G$ has
torsion index $1$, so we can use any coefficient ring $\k$.  The
Schubert polynomials are
\begin{gather*}
\lie{S}_{w_0}=\lie{S}=\eps_1^3\eps_2,\quad
\lie{S}_{s_2s_1s_2}=\eps_1^2\eps_2+\eps_1\eps_2^2,\quad
\lie{S}_{s_1s_2s_1}=\eps_1^3,\\
\lie{S}_{s_1s_2}=\eps_1^2+\eps_1\eps_2+\eps_2^2,\quad
\lie{S}_{s_2s_1}=\eps_1^2,\quad\lie{S}_{s_2}=\eps_1+\eps_2,\quad
\lie{S}_{s_1}=\eps_1,\quad\lie{S}_1=1,
\end{gather*}
where $s_1$ and $s_2$ denote the reflections in the roots
$\alpha_1=\eps_1-\eps_2$ and $\alpha_2=2\eps_2$.  The algebra $S^W$ is
freely generated by $p_1=\eps_1^2+\eps_2^2$ and
$p_2=\eps_1^2\eps_2^2$.  The matrices of $s_1$ and $s_2$ are
{\scriptsize$$
M_1=
\begin{pmatrix}
-1&0&0&0&0&0&0&0\\
0&1&-1&0&0&0&0&0\\
0&0&-1&0&0&0&0&0\\
p_1&0&0&1&0&0&0&0\\
0&0&0&0&-1&0&0&0\\
0&0&p_1&0&0&1&1&0\\
0&0&0&0&0&0&-1&0\\
-p_1^2&0&0&0&p_1&0&0&1
\end{pmatrix},\;
M_2=
\begin{pmatrix}
-1&0&0&0&0&0&0&0\\
0&-1&0&0&0&0&0&0\\
0&-2&1&0&0&0&0&0\\
0&0&0&-1&0&0&0&0\\
0&0&0&0&1&0&0&0\\
0&0&0&0&0&-1&0&0\\
0&2p_1&0&0&0&2&1&0\\
0&0&0&2p_1&0&0&0&1
\end{pmatrix},
$$
}and again one finds $A^W$ by solving the equations
$M_1\vec{c}=M_2\vec{c}=\vec{c}$.

\section{Homogeneous spaces}\label{section;homogeneous}

It is instructive to spell out some of the consequences of the results
of Section \ref{section;torsion} for the ordinary and equivariant
cohomology of homogeneous spaces.  As in the previous sections, $G$
denotes a compact connected Lie group with maximal torus $T$ and Weyl
group $W$, and $\k$ denotes a commutative ring with identity.  We
consider a fixed closed subgroup $U$ of $G$.

\begin{proposition}\label{proposition;equivariant-g/u}
Assume that $t(G)$ is a unit in $\k$.  Then there are $\k$-algebra
isomorphisms
$$
H_U^*(G/T;\k)\cong H_T^*(G/U;\k)\cong
S_\k\otimes_{(S^W)_\k}H^*(BU;\k).
$$
Hence $H_U^*(G/T;\k)$ is a free $H^*(BU;\k)$-module of rank $\abs{W}$
with generators of even degree.
\end{proposition}

\begin{proof}
Since $T$ and $U$ act freely on $G$, 
$$H_U^*(G/T;\k)\cong H_{U\times T}^*(G;\k)\cong H_T^*(G/U;\k).$$
The induction formula (see e.g.\ \cite[Section
1]{quillen;spectrum-equivariant-cohomology}) says that for any
$U$-space $Y$ there is a natural $\k$-algebra isomorphism
$H_U^*(Y;\k)\to H_G^*(G\times^UY;\k)$ which is linear over
$H^*(BG;\k)$.  Taking $Y=\pt$ gives $H^*(BU;\k)\cong H_G^*(G/U;\k)$.
We conclude, by Proposition
\ref{proposition;restriction}\eqref{item;tensor} and Corollary
\ref{corollary;bg-w}\eqref{item;bg-w}, that
$$
H_T^*(G/U;\k)\cong S_\k\otimes_{(S^W)_\k}H_G^*(G/U;\k)\cong
S_\k\otimes_{(S^W)_\k}H^*(BU;\k).
$$
The last assertion now follows from Proposition
\ref{proposition;restriction}\eqref{item;free}.
\end{proof}

Now let us assume $U$ to be connected and choose a maximal torus $T_U$
of $U$ which is contained in $T$.  We denote by $W_U$ the Weyl group
of $(U,T_U)$ and by $S_U$ the symmetric algebra of the character group
$\ca{X}(T_U)$.

\begin{corollary}\label{corollary;equivariant-g/u}
Assume that $U$ is connected and that both $t(G)$ and $t(U)$ are units
in $\k$.  Then $H_T^*(G/U;\k)$ is a free $\k$-module concentrated in
even degrees, finitely generated in each degree.  Moreover, there is a
$\k$-algebra isomorphism
$$
H_T^*(G/U;\k)\cong\bigl(S\otimes_{S^W}(S_U)^{W_U}\bigr)_\k.
$$
\end{corollary}

\begin{proof}
The first assertion follows from the fact that $H_T^*(G/U;\k)$ is free
and finitely generated over $H^*(BU;\k)$ with even degree generators
(Proposition \ref{proposition;equivariant-g/u}) and the fact that
$H^*(BU;\k)$ is free over $\k$ with even degree generators and
finitely generated in each degree (Corollary
\ref{corollary;bg-w}\eqref{item;bg-w}).  The second assertion follows
immediately from Proposition \ref{proposition;equivariant-g/u} by
applying Corollary \ref{corollary;bg-w}\eqref{item;bg-w} to the groups
$G$ and $U$, and by using
$$
(S\otimes\k)\otimes_{S^W\otimes\k}\bigl((S_U)^{W_U}\otimes\k\bigr)
\cong\bigl(S\otimes_{S^W}(S_U)^{W_U}\bigr)\otimes_\Z\k.
$$
\end{proof}

\begin{example}
Setting $U=T$ we obtain $H_T^*(G/T;\k)\cong(S\otimes_{S^W}S)_\k$ for
any $\k$ in which the torsion index of $G$ is invertible.
\end{example}

Next we record a few consequences of Theorem \ref{theorem;d}.

\begin{lemma}\label{lemma;g/u}
Assume that $U$ is connected and that $t(U)$ is a unit in $\k$.
\begin{enumerate}
\item\label{item;connected}
Let $I_U$ be the augmentation left ideal of the algebra of
$(S_U)^{W_U}\otimes\k$-linear endomorphisms of $S_U\otimes\k$.  Then
the projection $G/T_U\to G/U$ induces a $\k$-algebra isomorphism
$H^*(G/U;\k)\cong H^*(G/T_U;\k)^{I_U}$.
\item\label{item;maximal}
Assume that $U$ has maximal rank.  Let $\k_0=\Z[t(U)^{-1}]$.  Then
$$H^*(G/U;\k_0)\cong H^*(G/T;\k_0)^{W_U}.$$
Hence $H^*(G/U;\k)$ is a finitely generated free $\k$-module
concentrated in even degrees.
\end{enumerate}
\end{lemma}

\begin{proof}
Theorem \ref{theorem;d} implies that $H^*(X/U;\k)\cong
H^*(X/T_U;\k)^{I_U}$ for all free $U$-spaces $X$.  Taking $X=G$ proves
\eqref{item;connected}.  If $U$ has maximal rank, then $T_U=T$.  Since
the ring $\k_0=\Z[t(U)^{-1}]$ is principal of characteristic $0$, the
isomorphism in \eqref{item;maximal} follows from Proposition
\ref{proposition;kernel-cokernel}\eqref{item;cokernel} and the fact
that $H_T^*(G;\k_0)=H^*(G/T;\Z)\otimes_\Z\k_0$ is a finitely generated
free $\k_0$-module concentrated in even degrees.  Hence
$H^*(G/U;\k_0)$, being isomorphic to a submodule of $H^*(G/T;\k_0)$,
is likewise finitely generated, free and concentrated in even degrees.
Therefore, by the universal coefficient theorem, the same is true for
$H^*(G/U;\k)\cong H^*(G/U;\k_0)\otimes_{\k_0}\k$.
\end{proof}

This result is not optimal in all cases.  For instance, it is known
that if $U$ is the centralizer of a subtorus of a maximal torus $T$ of
$G$, then $H^*(G/U;\Z)\cong H^*(G/T;\Z)^{W_U}$, regardless of the
torsion index of $U$.  (A proof can be found in \cite[Theorem
III$''$]{bott-samelson;applications-theory-morse-symmetric} or
\cite[Theorem 5.5]{bernstein-gelfand-gelfand;schubert}.)

For subgroups of maximal rank we can now give a presentation of the
ordinary cohomology of the homogeneous space.  For coefficients in a
field $\k$, the isomorphism \eqref{equation;presentation} can be found
in Baum's paper \cite{baum;cohomology-homogeneous-spaces} and
Gugenheim and May's monograph
\cite{gugenheim-may;differential-torsion}.

\begin{theorem}\label{theorem;maximal-rank}
Assume that $U$ is connected and of maximal rank and that both $t(G)$
and $t(U)$ are units in $\k$.  Then we have $\k$-algebra isomorphisms
\begin{align}
\label{equation;equivariant-g/u}
H_T^*(G/U;\k)&\cong(S\otimes_{S^W}S^{W_U})_\k,\\
\label{equation;presentation}
H^*(G/U;\k)&\cong\bigl(S^{W_U}\big/S^W_+\bigr)_\k,
\end{align}
where $S^W_+$ denotes the ideal of $S^{W_U}$ generated by the
$W$-invariant elements of $S$ of positive degree.
\end{theorem}

\begin{proof}
The isomorphism \eqref{equation;equivariant-g/u} follows from
Corollary \ref{corollary;equivariant-g/u}.  For the isomorphism
\eqref{equation;presentation} we use $H^*(BG;\k)\cong(S^W)_\k$ and
$H^*(BU;\k)\cong(S^{W_U})_\k$ to get
$$
\bigl(S^{W_U}\big/S^W_+\bigr)_\k\cong(S^{W_U})_\k\otimes_{(S^W)_\k}\k
\cong H^*(BU;\k)\otimes_{H^*(BG;\k)}\k.
$$
It follows from Corollary \ref{corollary;bg-w}\eqref{item;bg-w} and
Lemma \ref{lemma;g/u}\eqref{item;maximal} that the cohomology groups
of $BG$ and of $G/U$ (with coefficients in $\k$) are free and vanish
in odd degrees.  This shows that the Leray spectral sequence of the
fibre bundle $G/U\to BU\to BG$ degenerates at the $E_2$ term, which
implies that $H^*(BU;\k)\otimes_{H^*(BG;\k)}\k$ is isomorphic to
$H^*(G/U;\k)$.
\end{proof}

(A comment on the hypotheses of Theorem \ref{theorem;maximal-rank}: in
all examples of maximal rank subgroups $U$ that we have checked, every
prime factor of $t(U)$ is also a prime factor of $t(G)$, and therefore
$t(U)$ is a unit in $\k$ if $t(G)$ is.  But we do not know if this is
true for all $U$.)

If $G=\SU(n,\C)$ or $G=\U(n,\H)$, then every closed connected maximal
rank subgroup $U$ has torsion index $1$ (this follows for example from
the information in \cite[Table~5.1]{mimura-toda;topology-lie-groups}),
so Theorem \ref{theorem;maximal-rank} applies to any $U$ and to any
coefficient ring $\k$.  If $G$ is any covering group of $\SO(n,\R)$
with $n\ge3$, then the torsion index of every $U$ is a power of $2$,
so Theorem \ref{theorem;maximal-rank} applies to any $U$ and to any
coefficient ring $\k$ which is an algebra over $\Z[\frac12]$.

It follows from Theorem \ref{theorem;maximal-rank} that the
homogeneous space $G/U$ is ``$T$-equivariantly formal'' over $\k$:
every cohomology class on $G/U$ extends to a $T$-equivariant class.

\begin{corollary}\label{corollary;maximal-rank-ideal}
The notation and the hypotheses are as in Theorem
{\rm\ref{theorem;maximal-rank}}.  There is a $\k$-algebra isomorphism
$H^*(G/U;\k)\cong H_T^*(G/U;\k)/(S_+)_\k$, where $(S_+)_\k$ denotes
the ideal of $H_T^*(G/U;\k)$ generated by polynomials of positive
degree.
\end{corollary}

\begin{proof}
It follows from Theorem \ref{theorem;maximal-rank} that
\begin{multline*}
H^*(G/U;\k)\cong(S^{W_U}\otimes_{S^W}\Z)_\k
\cong(S^{W_U}\otimes_{S^W}S\otimes_S\Z)_\k\\
\cong H_T^*(G/U;\k)\otimes_{S_\k}\k\cong H_T^*(G/U;\k)/(S_+)_\k.
\end{multline*}
\end{proof}

The fibre bundle $G/U\to BU\to BG$ is the right edge of the
commutative diagram
\begin{equation}\label{equation;four-bundles}
\vcenter{\xymatrix{
{\ar^j(0,0)*+{U/T};(12,16)*+{G/T}}
{\ar^q(12,16)*+{G/T};(24,32)*+{G/U}}
{\ar^{\bar{q}}(36,16)*+{BU};(48,0)*+{BG,}}
{\ar^{\bar{\jmath}}(24,32)*+{G/U};(36,16)*+{BU}}
{\ar^{i_U}(0,0)*+{U/T};(24,10.67)*+{BT}}
{\ar^{p_U}(24,10.67)*+{BT};(36,16)*+{BU}}
{\ar^i(12,16)*+{G/T};(24,10.67)*+{BT}}
{\ar^p(24,10.67)*+{BT};(48,0)*+{BG,}}
}}
\end{equation}
in which each of the four straight paths represents a locally trivial
fibre bundle.  As we have seen, under the hypotheses of Theorem
\ref{theorem;maximal-rank} the Leray spectral sequences of the two
bundles at the bottom and of the right edge collapse.  Our next result
gives a formula for a section of the restriction homomorphism
$\bar{\jmath}^*$ induced by the inclusion map $\bar{\jmath}\colon
G/U\inj BU$ indicated in \eqref{equation;four-bundles}.  In the course
of proving this formula we will see that the restriction map $j^*$
induced by the inclusion $j\colon U/T\to G/T$ is surjective, thus
proving that the spectral sequence of the left edge collapses as well.

\begin{lemma}\label{lemma;maximal-rank}
The notation and the hypotheses are as in Theorem
{\rm\ref{theorem;maximal-rank}} and in the diagram
\eqref{equation;four-bundles}.  A section of $\bar{\jmath}^*$ is given
by the map
$$
s_U=s\circ q^*\colon H^*(G/U;\k)\to(S^{W_U})_\k,
$$
where $s\colon H^*(G/T;\k)\to(S^W)_\k$ is the section of $i^*$ defined
in \eqref{equation;split}.  As an $(S^W)_\k$-module, $(S^{W_U})_\k$ is
free of rank $\abs{W}/\abs{W_U}$.
\end{lemma}

\begin{proof}
Let $\lie{S}_w\in H^*(BT;\k)$ be the classes defined in
\eqref{equation;classes}, but with $w\in W_U$ and with the operators
$\partial_w$ defined relative to a set of positive roots in the root
system of $(U,T)$.  Then the classes $(i_U^*(\lie{S}_w))_{w\in W_U}$
form a $\k$-basis of $H^*(U/T;\k)$ and, since $i\circ j=i_U$, they are
the restrictions of the classes $(i^*(\lie{S}_w))_{w\in W_U}$ to $U/T$
under the inclusion $j$.  Therefore $j^*$ is surjective.  Moreover,
$H^*(U/T;\k)$ is free and finitely generated, so the Leray-Hirsch
theorem applies to the left edge of \eqref{equation;four-bundles}.  We
conclude that there are isomorphisms
\begin{gather}
H^*(BT;\k)\cong H^*(BU;\k)\otimes_\k H^*(U/T;\k),\notag\\
\label{equation;g/u/t}
H^*(G/T;\k)\cong H^*(G/U;\k)\otimes_\k H^*(U/T;\k),
\end{gather}
expressing $H^*(BT;\k)$ as a free module over $H^*(BU;\k)$ with basis
$(\lie{S}_w)_{w\in W_U}$, and $H^*(G/T;\k)$ as a free module over
$H^*(G/U;\k)$ with basis $(i^*(\lie{S}_w))_{w\in W_U}$.  Now let $z\in
H^*(G/U;\k)$ and let $a=s_U(z)=s(q^*(z))\in H^*(BT;\k)$.  Then
$i^*(a)=q^*(z)$.  Expanding $a=\sum_{w\in W_U}x_w\lie{S}_w$ in terms
of the basis elements $\lie{S}_w$ with coefficients $x_w\in
H^*(BU;\k)$ and applying $i^*$ gives
\begin{equation}\label{equation;expansion}
q^*(z)=i^*(a)=\sum_{w\in W_U}i^*(x_w\lie{S}_w)=\sum_{w\in
W_U}\bar{\jmath}^*(x_w)i^*(\lie{S}_w).
\end{equation}
On the other hand, expanding $q^*(z)$ in terms of the basis elements
$i^*(\lie{S}_w)$ with coefficients in $H^*(G/U;\k)$ gives
$q^*(z)=z\mult1=z\mult i^*(\lie{S}_1)$.  Upon comparing this formula
with \eqref{equation;expansion} we find $\bar{\jmath}^*(x_w)=0$ for
$w\ne1$ and $\bar{\jmath}^*(x_1)=z$.  Hence $a=x_1$ is contained in
the submodule $H^*(BU;\k)$ of $ H^*(BT;\k)$ and $\bar{\jmath}^*(a)=z$.
This shows that $\bar{\jmath}^*s_U(z)=z$, i.e.\ $s_U$ is a section of
$\bar{\jmath}^*$.  Let $(a_i)_{i\in\ca{I}}$ be a basis of the
$\k$-module $H^*(G/U;\k)$, which is free and finitely generated by
Lemma \ref{lemma;g/u}\eqref{item;maximal}.  By \eqref{equation;g/u/t},
$H^*(G/T;\k)$ is free over $H^*(G/U;\k)$ of rank $\abs{W_U}$.
Moreover, $H^*(G/T;\k)$ is free over $\k$ of rank $\abs{W}$, so
$\abs{\ca{I}}=\abs{W}/\abs{W_U}$.  By the Leray-Hirsch theorem
(applied to the right edge of \eqref{equation;four-bundles}),
$(s_U(a_i))_{i\in\ca{I}}$ is a basis of the $(S^W)_\k$-module
$(S^{W_U})_\k$, and therefore $(S^{W_U})_\k$ is free of rank
$\abs{W}/\abs{W_U}$.
\end{proof}

We can now prove versions of Proposition \ref{proposition;restriction}
and Theorem \ref{theorem;d} relative to the subgroup $U$:
$G$-equivariant cohomology determines $U$-equivariant cohomology and
vice versa.  As before, $X$ denotes a topological $G$-space.

\begin{theorem}\label{theorem;u}
The hypotheses are as in Theorem {\rm\ref{theorem;maximal-rank}}.  Let
$I_{G,U}$ be the augmentation left ideal of the algebra of
$(S^W)_\k$-linear endomorphisms of $(S^{W_U})_\k$.  Then we have
$\k$-algebra isomorphisms
\begin{gather}
\label{equation;u-equivariant}
H_U^*(X;\k)\cong(S^{W_U})_\k\otimes_{(S^W)_\k}H_G^*(X;\k),\\
\label{equation;u-equivariant-delta}
H_G^*(X;\k)\cong H_U^*(X;\k)^{I_{G,U}}.
\end{gather}
\end{theorem}

\begin{proof}
Define
$$\bar{s}_U\colon H^*(G/U;\k)\to H_U^*(X;\k)$$
by $\bar{s}_U=\pr_U^*\circ s_U$, where $\pr_U\colon X_U\to BU$ is the
natural projection and $s_U$ is as in Lemma \ref{lemma;maximal-rank}.
Then $\bar{s}_U$ is a section of the restriction map $H_U^*(X;\k)\to
H^*(G/U;\k)$.  Let $\{a_i\}_{i\in\ca{I}}$ be a $\k$-basis of
$H^*(G/U;\k)$ as in the proof of Lemma \ref{lemma;maximal-rank}.
Applying the Leray-Hirsch theorem to the fibre bundle $G/U\to X_T\to
X_U$, we conclude that $\{\bar{s}_U(a_i)\}_{i\in\ca{I}}$ is a basis of
$H_U^*(X;\k)$ as a $H_G^*(X;\k)$-module.  The natural homomorphism
$$(S^{W_U})_\k\otimes_{(S^W)_\k}H_G^*(X;\k)\longto H_U^*(X;\k)$$
maps the basis element $s_U(a_i)\otimes1$ of the module on the left to
the basis element $\bar{s}_U(a_i)$ of the module on the right and
therefore is an isomorphism.  This proves
\eqref{equation;u-equivariant}.  It follows from Lemma
\ref{lemma;maximal-rank} that $(S^{W_U})_\k$ is a progenerator of the
category of $(S^W)_\k$-modules.  The isomorphism
\eqref{equation;u-equivariant-delta} now follows from
\eqref{equation;u-equivariant} by applying the first Morita
equivalence theorem as in the proof of Theorem \ref{theorem;d}.
\end{proof}

In particular, we find a $U$-equivariant version of the isomorphism
\eqref{equation;equivariant-g/u}.

\begin{corollary}\label{corollary;u}
The hypotheses are as in Theorem {\rm\ref{theorem;maximal-rank}}.  We
have a $\k$-algebra isomorphism
$$
H_U^*(G/U;\k)\cong(S^{W_U}\otimes_{S^W}S^{W_U})_\k.
$$
\end{corollary}

\begin{proof}
This follows immediately from Theorem \ref{theorem;u} by taking
$X=G/U$ and using $H_G^*(G/U;\k)\cong H^*(BU;\k)\cong(S^{W_U})_\k$.
\end{proof}

\begin{corollary}\label{corollary;u-root}
The hypotheses are as in Theorem {\rm\ref{theorem;maximal-rank}}.  In
addition, if $G$ has a root $\alpha$ such that $\alpha/2\in\X(T)$,
assume that $2$ is not a zero divisor in $\k$.  Then we have
$\k$-algebra isomorphisms
\begin{align*}
H_U^*(G/U;\k)&\cong(S_\k)^{W_U}\otimes_{(S_\k)^W}(S_\k)^{W_U},\\
H_T^*(G/U;\k)&\cong S_\k\otimes_{(S_\k)^W}(S_\k)^{W_U},\\
H^*(G/U;\k)&\cong(S_\k)^{W_U}\big/I_\k.
\end{align*}
\end{corollary}

\begin{proof}
The root system of $(U,T)$ is a subset of the root system of $(G,T)$,
so Corollary \ref{corollary;abelian} applies to both $G$ and $U$,
telling us that $(S^W)_\k=(S_\k)^W$ and $(S^{W_U})_\k=(S_\k)^{W_U}$.
The result now follows from Theorem \ref{theorem;maximal-rank} and
Corollary \ref{corollary;u}.
\end{proof}

\begin{example}
Let $G=\U(n,\H)$ and let $U$ be the diagonal subgroup $\U(1,\H)^n$.
Then $t(G)=t(U)=1$, so Corollary \ref{corollary;u} applies to any
coefficient ring $\k$.  The space $G/U$ is the variety of full flags
of the quaternionic vector space $\H^n$.  Let $T=\U(1,\C)^n$ be the
diagonal maximal torus of $G$ and $U$.  Then $W_U\cong(\Z/2\Z)^n$ acts
by sign changes on each of the coordinates of $\t\cong\R^n$, and
$W\cong W'\ltimes W_U$ acts by signed permutations on $\t$, where
$W'=S_n$.  We have $S=\Z[y_1,y_2,\dots,y_n]$, a polynomial ring in
variables $y_i$ of degree $2$.  Let $S'=S^{W_U}$ and $x_i=y_i^2$.
Then $S'=\Z[x_1,x_2,\dots,x_n]$, so we obtain
$$
H_U^*(G/U;\Z)\cong S'\otimes_{(S')^{W'}}S'\cong(S'\otimes_\Z S')/I.
$$
Here $I$ is the ideal of $S'\otimes_\Z S'$ generated by all elements
of the form $u\otimes1-1\otimes u$, where $u\in S'$ is a symmetric
polynomial in the variables $x_i$.  Cf.\ Mare's paper
\cite{mare;equivariant-quaternionic-flag}, where this algebra is
computed by a different method.  The cohomology $H_U^*(G/U;\k)$ is now
found by tensoring $(S'\otimes_\Z S')/I$ with $\k$.  Observe also that
$H_U^*(G/U;\k)$ is \emph{not} isomorphic to
$(S_\k)^{W_U}\otimes_{(S_\k)^W}(S_\k)^{W_U}$ if $\k$ has
characteristic $2$.  This shows that the conditions of Corollary
\ref{corollary;u-root} cannot be weakened.
\end{example}


\bibliographystyle{amsplain}

\def\cprime{$'$}
\providecommand{\bysame}{\leavevmode\hbox to3em{\hrulefill}\thinspace}
\providecommand{\MR}{\relax\ifhmode\unskip\space\fi MR }
\providecommand{\MRhref}[2]{%
  \href{http://www.ams.org/mathscinet-getitem?mr=#1}{#2}
}
\providecommand{\href}[2]{#2}


\end{document}